\documentstyle{amsppt}
\magnification=1200

\pagewidth{15 truecm}
\pageheight{23 truecm}
\baselineskip=12pt
\parindent=20pt
\TagsOnRight
\def\spb{\smallpagebreak}
\def\mpb{\vskip 0.5truecm}
\def\bpb{\vskip 1truecm}
\hcorrection{.25in}
\NoBlackBoxes

\topmatter

\title
Rigidity of Irreducible Hermitian Symmetric Spaces
of the Compact Type under K\"ahler Deformation
\endtitle

\rightheadtext{Rigidity of Certain Symmetric Spaces
under K\"ahler Deformation\quad}

\author
Jun-Muk Hwang\\
Ngaiming Mok
\endauthor

\thanks
Hwang was partially supported by GARC-KOSEF, and at
MSRI by NSF grant DMS 9022140.
Mok was partially supported by a grant 
from The University of Hong Kong.
\endthanks

\endtopmatter
\document

\abstract
We study deformations of irreducible Hermitian symmetric spaces $S$
of the compact type, known to be locally rigid, as projective-algberaic
manifolds and prove that no jump of complex structures can occur.
For each $S$ of rank $\ge 2$ there is an associated reductive linear
group $G$ such that $S$ admits a holomorphic $G$-structure,
corresponding to a reduction of the structure group of the tangent bundle.
$S$ is characterized as the unique simply-connected compact complex
manifold admitting such a $G$-structure which is at the same time
integrable. To prove the deformation rigidity of $S$ it suffices that
the corresponding integrable $G$-structures converge.

We argue by contradiction using the deformation theory of rational curves.
Assuming that a jump of complex structures occurs, cones of vectors
tangent to degree-1 rational curves on the special fiber $X_0$ are linearly
degenerate, thus defining a proper meromorphic distribution $W$ on $X_0$.
We prove that such $W$ cannot possibly exist. On the one hand,
integrability of $W$ would contradict the fact that $b_2(X)=1$.
On the other hand, we prove that $W$ would be automatically integrable
by producing families of integral complex surfaces of $W$ as pencils
of degree-1 rational curves. For the verification that there are
enough integral surfaces we need a description of generic 
cones on the special fiber. We show that they are in fact images of
standard cones under linear projections. We achieve this by studying
deformations of normalizations of Chow spaces of minimal rational curves
marked at a point, which are themselves Hermitian symmetric, irreducible
except in the case of Grassmannians.
\endabstract

Let $S$ be a Hermitian symmetric space of the compact type. It follows from
Bott [Bo] that the complex structure of $S$
is infinitesimally rigid. In the special  case where $S=\Bbb P^1\times\Bbb P^1$
is a product of two Riemann spheres, it is well-known that $S$ can be
holomorphically deformed to any Hirzebruch surface $\Sigma_a$ with $a$ even.
In the case where $S$ is irreducible, following Kodaira-Spencer [KS, 1958]
it was however conjectured that $S$ is rigid under deformation as a 
complex manifold. While this was affirmed in the case of the projective 
space $\Bbb P^n$ by Siu [S1, 1992] and in the case of the hyperquadric $Q^n$, 
$n\ge 3$, by Hwang [H1, 1995], the general case remains
open. In the present article we consider the special case of K\"ahler
deformations and prove the rigidity of $S$ under such deformations. More
precisely, we prove

\spb
\proclaim{Theorem 1} Let $S$ be an irreducible Hermitian symmetric space of
the compact type. Let $\pi:X\to\triangle$ be a
regular family of compact complex manifolds over the unit disk $\triangle$. 
Suppose $X_t:=\pi^{-1}(t)$ is biholomorphic to $S$ for $t\ne 0$ 
and the central fiber $X_0$ is K\"ahler. 
Then, $X_0$ is also biholomorphic to $S$.
\endproclaim

\spb
Since the second Betti number of $S$ is 1, the K\"ahler condition on $X_0$ is 
\linebreak
equivalent to the assumption that $X_0$ is projective-algebraic.
We have the following equivalent algebro-geometric formulation of our result
(cf. (5.2)).

\spb
\proclaim{Theorem 1$'$} Let $S$ be an irreducible Hermitian symmetric space of
the compact type. Let $\rho:X\to Z$ be a smooth proper morphism between
two connected algebraic varieties over $\Bbb C$. Suppose for some 
point $y$ on $Z$, the fiber $X_y$ is biregular to $S$.
Then, for any point $z$ on $Z$, the fiber $X_z$ is biregular to $S$.
\endproclaim

\spb
Let $M$ be a compact K\"ahler manifold homeomorphic to the projective
space $\Bbb P^n$ (resp. the hyperquadric $Q^n$, $n\ge 3$) such that
the canonical line bundle is not ample. Then, by the work of 
Hirzebruch-Kodaira [HK, 1957] (resp. Brieskorn [B, 1964]) using the
Riemann-Roch formula, it was well-known that $M$ is biholomorphic to
$\Bbb P^n$ (resp. to $Q^n$, $n\ge 3$). Theorem 1 for $\Bbb P^n$ and $Q^n$
$(n\ge 3)$ already follows from these early results. The difficulty of
our Theorem 1 lies therefore in the case of other model spaces $S$,
for which there do not exist algebro-geometric characterizations in terms
of holomorphic line bundles. As such our difficulties are quite distinct
from those encountered in Siu [S1] and Hwang [H1, 2], where the major obstacle
was the lack of the K\"ahler condition on the central fiber.

\spb
For a historical overview on the general problem of rigidity
(or nonrigidity) of complex structures of Hermitian symmetric
spaces of the compact type, and more specifically on the problem of
deformation (non)rigidity of such spaces, we refer the reader to Siu [S2].
Regarding our Theorem 1, Mabuchi [Ma] kindly provided the
first author with a preprint of his, explaining a proof using Euler vector
fields. As far as we are aware, it has so far not been possible to complete
his scheme of proof except in some special cases such as Grassmannians of
$\text{rank}\,2$.

\spb
Our method of proof is based on two ingredients: the analytic ingredient
of holomorphic $G$-structures and the algebro-geometric ingredient of
deformation of rational curves. Hermitian symmetric spaces of the compact
type $S$ are compactifications of Euclidean spaces, obtained canonically 
using the Harish-Chandra Embedding Theorem. 
For $S$ irreducible and of $\text{rank}\ge 2$, 
$S$ admits an integrable (holomorphic) $K^{\Bbb C}$-structure
for some proper complex Lie subgroup $K^{\Bbb C}$ of the general linear group.
Here a $K^{\Bbb C}$-structure refers to a holomorphic reduction
of the holomorphic frame bundle from the general linear group to $K^{\Bbb C}$.
Integrability refers to the possibility of realizing such a
reduction by means of holomorphic local coordinates on $S$. 
By a result of Ochiai [Oc], 
a simply-connected compact complex manifold with an integrable
$K^{\Bbb C}$-structure is biholomorphic to the model space $S$.
This may be taken as our point of departure, reducing our Theorem 1 to
proving the convergence of the integrable $K^{\Bbb C}$-structures on
$X_t$, $t\ne 0$, to one on the central fiber $X_0$. As integrability is a
closed condition, the problem is reduced to verifying the
convergence of $K^{\Bbb C}$-structures.

\spb
Here enters the deformation theory of rational curves. On the model space
$S$ there exists at every $s\in S$ a homogeneous complex submanifold
$\Cal C_s\subset\Bbb P(T_s(S))$, called the cone at $s$, consisting of
all directions tangent to minimal rational curves. The collection of
$\Cal C_s$ then gives a holomorphic fiber bundle $\Cal C\mapsto S$
whose fibers are isomorphic to a model $\Cal C_o$. 
The holomorphic $K^{\Bbb C}$-structure
of $S$ can be recovered from the bundle of cones, $K^{\Bbb C}$ being
isomorphic to the group of linear transformations inducing an
isomorphism on $\Cal C_o$. We are thus led to examine the limit of
the bundles of cones over $X_t$, as $t\mapsto 0$.
To prove Theorem 1 it suffices to show that the limit
is a bundle of cones modeled on $\Cal C_o$
such that for $x\in X_0$, the embedding
$\Cal C_x\subset\Bbb P T_x(X_0)$ is isomorphic to the
standard embedding $\Cal C_o\subset\Bbb P T_o(S)$.

\spb
Remarkably the standard cone $\Cal C_o\subset\Bbb P T_0(S)$
is itself a Hermitian symmetric space of the
compact type, of $\text{rank}\le 2$, irreducible except in the case of
Grassmannians of $\text{rank}\ge 2$. 
This entitles us to an argument by induction. 
For the model space $S$, $\Cal C_s$ is isomorphic to the Chow 
space $\Cal M_s$ of minimal rational curves marked at $s$. From
the K\"ahler condition on $X_0$ we deduce that the limiting normalized 
Chow space $\Cal M_x$ is smooth at a generic point $x\in X_0$. For $S$ 
different from Grassmannians, $\Cal M_x$ is by induction a Hermitian symmetric
space. Embedding  $\Cal M_x\hookrightarrow\Bbb P(V_x)$ projectively, there is a
linear map $V_x\to T_x(X_0)$ inducing a rational map 
$\Phi_x:\Cal M_x\to \Cal C_x$. Our main task is reduced to proving that
$\Phi_x$ is an isomorphism. The same reduction is valid in the remaining
case of Grassmannians, for which the cones are products of two projective
spaces, as follows. From cohomological considerations we can exclude the 
possibility such as the deformation of $\Bbb P^1\times\Bbb P^1$ 
into non-trivial Hirzebruch surfaces. The point is to prove that
limits of direct factors of the cone cannot decompose. As deformations
of such direct factors arise from deformations of projective
subspaces of Grassmannians, we are led to proving that limits of
such projective spaces cannot decompose, which follows from the
cohomological structure of Grassmannians and the K\"ahler condition
on~$X_0$.

\spb
The main difficulty in proving Theorem 1 is to rule out the possibility that
the cones $\Cal C_x\subset\Bbb P T_x(X_0)$ are linearly 
degenerate. In this case we have a meromorphic
distribution $W$ on $X_0$, defined at generic points as the linear
span of vectors tangent to minimal rational curves. If this distribution
is integrable where holomorphic, we show that leaves of the
associated foliation can be compactified to
subvarieties of $X_0$, and obtain a contradiction to the fact that
the Picard number of $X_0$ is 1, using the deformation theory of rational
curves.

\spb
Our theorem is therefore finally reduced to another problem of integrability,
one for which we need to verify the Frobenius condition for the distribution
$W$ defined algebro-geometrically using minimal rational curves. 
We further reduce the verification to finding at a generic point $x$ of $X_0$
a certain set of integral surfaces $\Sigma$ such that 
$\Lambda^2 T_x(\Sigma)$ generates $\Lambda^2 W_x$. Such integral surfaces
can be obtained as pencils of minimal rational curves,
and the linear-algebra statement can be verified using a uniform
description of cones as Zariski closures of graphs of quadratic
maps on Euclidean spaces. Alternatively it can also be verified
using descriptions of second exterior powers of isotropy representations 
(cf. Onishchik-Vinberg [OV]).

\spb
In \S 1 we collect relevant basic facts about $G$-structures
arising from irreducible Hermitian symmetric spaces of $\text{rank}\ge 2$,
including a discussion on integrability. We also prove a Hartogs-type
extension result on holomorphic $G$-structures for $G$ reductive,
allowing us in what follows to concentrate our discussion at generic
points of the central fiber. In \S 2 we describe the cones of minimal
rational curves on the model spaces $S$, noting that a model cone
$\Cal C_o$ is itself Hermitian symmetric, irreducible except in the
case of Grassmannians. We observe the generic smoothness of normalized 
Chow spaces of marked minimal rational curves on the central fiber $X_0$, 
and prove that generic minimal rational curves on $X_0$ 
are immersed and of the standard type 
(i.e., direct factors of normal bundles are of degree 1 or 0).
In \S 3 we prove that in the case of deformations of Grassmannians,
generic limits of normalized Chow spaces of marked minimal rational curves are
biholomorphic to the standard cone. 
In \S 4 we consider distributions $W$ on Fano manifolds spanned at
generic points by tangents to standard minimal rational curves.
We prove in particular that $W$ is integrable whenever the varieties
of tangent lines to cones are linearly non-degenerate in
$\Bbb P\Lambda^2 W$ at generic points. In \S 5 we apply results of \S 3 and
\S 4 to our central fiber $X_0$ to deduce the linear non-degeneracy
of cones at generic points, completing the proofs of Theorems 1 and 1$'$.

\newpage

\bpb\noindent
\underbar{Table of Contents}

\spb
\item{\S 1} $G$-structures associated to Hermitian symmetric spaces

\spb
\item{\S 2} Chow spaces and cones of minimal rational curves

\spb
\item{\S 3} Non-deformability of normalized 
Chow spaces in the case of Grassmannians

\spb
\item{\S 4} The linear span of tangents to minimal rational curves
in Fano maifolds

\spb
\item{\S 5} Cones in the central fiber

\bpb\noindent
{\bf \S 1 $G$-structures associated to Hermitian symmetric spaces}

\mpb\noindent
(1.1)\quad Let $n$ be a positive integer.
Fix an $n$-dimensional complex vector space $V$ and let
$M$ be any $n$-dimensional complex manifold. In what follows all
bundles are understood to be holomorphic. The frame bundle $\Cal F(M)$
is a principal $GL(V)$-bundle with the fiber at $x$ defined as
$\Cal F(M)_x=\text{Isom}(V,T_x(M))$, the set of linear isomorphisms
from $V$ to the holomorphic tangent space at $x$. Let $G\subset GL(V)$
be any complex Lie subgroup. A (holomorphic) $G$-structure is a $G$-principal
subbundle $\Cal G(M)$ of $\Cal F(M)$. For $G\ne GL(V)$ we say that
$\Cal G(M)$ defines a holomorphic reduction of the tangent bundle to $G$.

\spb
Let $\varphi_\alpha:U_\alpha\to V$ be a chart on $M$. In terms of
Euclidean coordinates we identify $\Cal F(U_\alpha)$ with the product
$GL(V)\times U_\alpha$. We say that a $G$-structure
$\Cal G(M)$ on $M$ is integrable if and only if there exists an atlas of
charts $\{\varphi_\alpha:U_\alpha\to V\}$ such that the restriction
$\Cal G(U_\alpha)$ of $\Cal G(M)$ to $U_\alpha$ is the product 
$G\times U_\alpha\subset GL(V)\times U_\alpha$. 

\mpb\noindent
(1.2)\quad Let $(S,h)$ be a Hermitian symmetric space of the compact type
and $G$ be the identity component of the isometry group of $(S,h)$.
(The meaning of $G$ here is not to be confused with that appearing
in ``$G$-structure'', where $G$ is used as a generic symbol for groups.)
Let $G^{\Bbb C}$ be the identity component of the group of automorphisms
of $S$. Then, $G^{\Bbb C}$ is a semisimple complex Lie group and
$G\subset G^{\Bbb C}$ is a compact real form. Let $K\subset G$ be
the isotropy subgroup at a reference 
point $o\in S$, and $\frak k\subset\frak g$ be corresponding Lie algebras.
Write $\frak g=\frak k+\frak m$ for the Cartan decomposition, where
$\frak m$ is the real tangent space at $o$. For the complexifications
$\frak g^{\Bbb C}=\frak k^{\Bbb C}+\frak m^- + \frak m^+$, where
$\frak m^{\Bbb C}=\frak m^- + \frak m^+$ is the decomposition of
the complexified tangent space into the direct sum of subspaces 
$\frak m^-$ resp. $\frak m^+$ of type (0, 1) resp. (1, 0). Let
$P\subset G^{\Bbb C}$ be the isotropy subgroup at $o$. We have
$\frak p = \frak k^{\Bbb C} + \frak m^-$ at the level of Lie algebras.
$\frak m^+$, $\frak m^-\subset \frak g^{\Bbb C}$ are abelian subalgebras.
The corresponding Lie subgroups $M^+$, $M^-\subset G^{\Bbb C}$ are
isomorphic to the abelian Lie group $\Bbb C^n$. According to the
Harish-Chandra Embedding Theorem, the map $M^+\times K^{\Bbb C}\times M^-
\to G^{\Bbb C}$ defined by $(m^+,k,m^-)\to m^+\cdot k\cdot m^-$ is
an open embedding. Applying to $o$ we have an embedding of the
$M^+$-orbit of $o$ into $S$. Identifying $M^+$ with $\Bbb C^n$, taking $o$
to be the origin, we regard $S$ as a compactification of $\Bbb C^n$ such
that translations on $\Bbb C^n$ extend to automorphisms of $S$.
Euclidean coordinates on $\Bbb C^n\subset S$ will be referred to as 
Harish-Chandra coordinates.

\spb
On $\Bbb C^n$, $K^{\Bbb C}$ acts as a group of linear transformations,
acting faithfully on $T_o(S)$. Denote the corresponding subgroup of
$GL(T_o(S))$ by $K_o^{\Bbb C}$. For $s\in S$ denote by $P_s\subset G^{\Bbb C}$
the isotropy subgroup at $s$. Consider the homomorphism
$\varphi_s:P_s \to GL(T_s(S))$ defined by $\varphi_s(f)=df(s)$.
At $o$, $\text{Ker}(\varphi_o)=M^-$ and $\varphi_o(P_o)=K_o^{\Bbb C}$.
For $x\in\Bbb C^n\subset S$ we use the trivialization of $T(\Bbb C^n)$
via Harish-Chandra coordinates to identify $T_x(S)$ with $T_o(S)\cong\Bbb C^n$.
Since $P_s$ is conjugate to $P_o=P$ via a Euclidean translation in
$M^+$, all $K_x^{\Bbb C}$ are identical under the corresponding 
identification $GL(T_x(S))\cong GL(n,\Bbb C)$. Using a finite number of
Harish-Chandra coordinate charts, we obtain on $S$ an integrable 
$G$-structure with $G=K_o^{\Bbb C}$. From now on we will identify
$K_o^{\Bbb C}$ with $K^{\Bbb C}$ and call this the $K^{\Bbb C}$-structure
on $S$. For $S$ of $\text{rank}\ge 2$, $K^{\Bbb C}\ne GL(V)$ and
we have a non-trivial reduction of the holomorphic tangent bundle.

\spb
For details on the Harish-Chandra Embedding Theorem, we refer the reader
to Wolf [W].

\mpb\noindent
(1.3)\quad Let $S$ be an irreducible Hermitian symmetric manifold
of the compact type and of $\text{rank}\ge 2$. For the regular family
$\pi:X\to\triangle$ as in Theorem 1 we will be considering limits
of $K^{\Bbb C}$-structures on $X_t$, $t\ne 0$, as $t\to 0$. In order to
have a limit defined everywhere on $X_0$, we will need the following
Hartogs-type extension result. By a holomorphic family $\sigma:M\to\triangle$
of complex manifolds $M_t$ we mean a holomorphic submersion $\sigma$ with
fibers $M_t:=\sigma^{-1}(t)$. Denote by $\Cal F_\sigma(M)$ the
relative frame bundle where 
$\Cal F_{\sigma,x}(M)=\text{Isom}(V,T_{\sigma,x}(M))$,
$T_{\sigma,x}(M)=\text{Ker}(d\sigma(x))$. For a subgroup 
$G\subset GL(V)$, by a holomorphic family of $G$-structures on $M$
we mean a holomorphic $G$-principal subbundle $\Cal G$ of $\Cal F_\sigma(M)$.
The restriction to $M_t$ will be denoted by $\Cal G_t$. We have

\spb
\proclaim{Proposition 1} Let $H\subset GL(V)$ be a reductive algebraic
subgroup, $\sigma:M\to\triangle$ be a holomorphic family of 
$n$-dimensional complex manifolds $M_t$. Let $E\subset M_0$ be a
complex-analytic subvariety of codimension $\ge 1$ and $\Cal H\subset 
\Cal F_\sigma(M-E)$ be a holomorphic family of $H$-structures. Then, $\Cal H$
extends holomorphically to $M$.
\endproclaim

\spb
\demo{Proof} The problem being local we may assume 
$M=\triangle^n\times\triangle$ with $\sigma(x,t)=t$.
$\Cal H\subset\Cal F_\sigma((\triangle^n\times\triangle)-E)$ is
equivalently given by a holomorphic map $f:(\triangle^n\times\triangle)-E\to
GL(V)/H$ into the complex homogeneous space $GL(V)/H$, 
$E\subset\triangle^n\times\triangle$ being of complex codimension $\ge 2$.
By Matsushima-Morimoto [MM], $GL(V)/H$ is Stein for $H$ reductive.
Embedding $GL(V)/H$ as a complex submanifold of some Euclidean space and
applying Hartogs' Extension Theorem for holomorphic functions to
$(\triangle^n\times\triangle)-E$, we conclude that $f$ and hence $\Cal H$
extend holomorphically to $\triangle^n\times\triangle$, as desired.
\enddemo

\mpb\noindent
(1.4)\quad For the proof of Theorem 1, once a limit holomorphic
$K^{\Bbb C}$-structure is constructed on the central fiber it will be
immediate that the $K^{\Bbb C}$-structure is integrable, by the closedness
of the integrability condition. To give simple
references, we sketch a proof in the special case of
$G$-structures of finite type. To each (holomorphic) $G$-structure 
$\Cal G\subset \Cal F(M)$ on a complex manifold $M$ we can associate
a prolongation bundle $\Cal P_1\to\Cal G$, which is a $G_1$-structure
on $\Cal G$ for some $G_1\subset GL(V\oplus\frak g)$, $\frak g$ denoting
the Lie algebra of $G$ (cf. Sternberg [St]). The Lie algebra
$\frak g_1\subset \frak g\frak l(V\oplus \frak g)$, called the first
prolongation algebra, is constructed from the embedding
$\frak g\subset\frak g\frak l (V)$. Higher prolongation algebras
$\frak g_{k+1}=(\frak g_k)_1$ and prolongation bundles
$\Cal P_{k+1}\to\Cal P_k$ can be constructed inductively. The
$G$-structure is said to be of finite type if $\frak g_m=0$ for some $m$,
of type $k$ if $k$ is the first index for which $\frak g_k=0$. 
$K^{\Bbb C}$-structure associated to irreducible Hermitian symmetric
spaces $S$ of $\text{rank}\ge 2$ are of type 2 (cf. Ochiai [Oc]). We have

\spb
\proclaim{Proposition 2} Suppose that we have a holomorphic family $\Cal G$
of $G$-structures of finite type on a family $\sigma:M\to\triangle$ of
complex manifolds $M_t$. Suppose that $\Cal G_t$ is integrable for all 
$t\ne 0$. Then, $\Cal G_0$ is also integrable.
\endproclaim

\spb
\demo{Proof} Under the finite type condition, by going to a prolongation
bundle we can assume that $G=\{e\}$ (cf. Sternberg
[St, p.338]). Thus, $\Cal G_t$ is a field of
frames on $M_t$, and it is integrable if the frame field corresponds
to a coordinate frame field. Since the problem is local, we may assume
that $M=\triangle^n\times\triangle$. By assumption, we have holomorphic
vector fields $V_1,\ldots,V_n$ on $M$, so that they define the
given frame field and $[V_i,V_j]=0$ for $t\ne 0$. But then they must
satisfy $[V_i,V_j]=0$ for $t=0$, too. Hence their integrals give a
coordinate system integrating the $\{e\}$-structure $\Cal G_0$.
\enddemo

\spb\noindent
{\smc Remarks}

\noindent
The integrability of $K^{\Bbb C}$-structures holds in a more general
setting than Proposition 2. In a forthcoming article, we will
show that a $K^{\Bbb C}$-structure on a Fano manifold is always integrable.

\spb
Finally, to recover Hermitian symmetric spaces from integrable 
$K^{\Bbb C}$-structures, we have

\spb
\proclaim{Proposition 3 (Ochiai [Oc])} 
Let $S$ be an irreducible Hermitian symmetric
space of the compact type and of $\text{rank}\ge 2$ with associated
$K^{\Bbb C}$-structures, $K^{\Bbb C}\subset GL(T_o(S))$; 
$o\in S$. Let $M$ be a compact simply-connected complex
manifold with an integrable $K^{\Bbb C}$-structure. Then, $M$ 
is biholomorphic to $S$.
\endproclaim

\bpb\noindent
{\bf \S 2 Chow spaces and cones of minimal rational curves}

\mpb\noindent
(2.1)\quad Let $S$ be an irreducible Hermitian symmetric space of
the compact type. The second homology group and hence the Picard
group of $S$ are infinite cyclic. Degrees of algebraic curves of
$S$ will be measured with respect to an ample line bundle $\Cal O(1)$
which is a generator of the Picard group. An algebraic curve
$C\subset S$ is called a minimal rational curve if and only if it is
of degree 1. One way of describing minimal rational curves is by means
of the first canonical embedding of $S$ into a projective space 
$\Bbb P^N$ (cf. Nakagawa-Takagi [NT]). Complex-analytically this is
nothing other than the embedding defined by $\Cal O(1)$, such as
the Pl\"ucker embedding of a Grassmannian. Identifying $S$ as a
subvariety of $\Bbb P^N$, the minimal rational curves on
$S$ are equivalently rational lines on $\Bbb P^N$ lying in $S$.
In particular, they are smooth. One can check the existence of such
rational lines using the definition of the first canonical embedding
via an embedding of complex Lie groups. The Chow space of all minimal
rational curves $C$ on $S$ constitutes a homogeneous space under the
action of $G^{\Bbb C}$, the identity component of the group 
$\text{\it Aut}(S)$ of automorphisms of $S$. For $s\in S$,  consider
the set $\Cal M_s$ of all minimal rational curves $C$ passing through
$s$ and the cone $\Cal C_s\subset\Bbb PT_s(S)$ defined as
$\{[\alpha]\in\Bbb P T_s(S): T_s(C)=\Bbb C\alpha$ for some
$[C]\in\Cal M_s\}$. The tangent map
$\Phi_s:\Cal M_s\to\Cal C_s$ thus defined is a bijection. We fix a
reference point $o\in S$ and call $\Cal C_o\subset\Bbb PT_o(S)$ the
standard cone associated to $S$. A remarkable fact about standard cones 
is that they are themselves Hermitian symmetric spaces of the compact type,
of $\text{rank}\le 2$, and irreducible except in the case of Grassmannians
of $\text{rank}\ge 2$, where $\Cal C_o$ is a product of two projective
spaces. This can be seen from the following table giving a classification
of $S$ and a tabulation of their standard cones, together with a 
description of the inclusion $\Cal C_o\subset \Bbb PT_o(S)$ as a projective
embedding of the Hermitian symmetric manifold $\Cal C_o$.
Here in this table alone $G$ will denote a finite simply-connected covering
of $\text{\it Aut}_0(S,h)$ and $K\subset G$ the isotropy subgroup.
$\Bbb O$ will stand for the octonions (Cayley numbers).

\mpb
\centerline{\underbar{Table of irreducible Hermitian symmetric spaces $S$}}
\centerline{\underbar{of the compact type and their standard cones $\Cal C_o$}}
$$
\vbox{\offinterlineskip\tabskip=0pt plus1fil
\halign to \displaywidth{\tabskip=0pt\strut
\hfil$#$\hfil\quad 
&\vrule#\
&\hfil$#$\hfil\
&\vrule#\
&\hfil$#$\hfil\
&\vrule#\
&\hfil$#$\hfil\
&\vrule#\ 
&\hfil$#$\hfil\
&\vrule#\
&\hfil$#$\hfil\
&\vrule#\
&\hfil$#$\hfil\
&\vrule#\
&$#$&\hfil$#$\hfil 
\tabskip=0pt plus1fil &\llap{$#$}\tabskip=0pt\cr 
\omit&\multispan{12}\hrulefill\cr
\omit&height3pt&\omit&&\omit&&\omit&&\omit&&\omit&&\omit&&\cr
&&\text{Type}&&G&&K&&G/K=S&&\Cal C_o&&\text{Embedding}&&&&\cr
\omit&height3pt&\omit&&\omit&&\omit&&\omit&&\omit&&\omit&&\cr
\omit&\multispan{12}\hrulefill\cr
\omit&height3pt&\omit&&\omit&&\omit&&\omit&&\omit&&\omit&&\cr
&&\text{I}&&SU(p+q)&&S(U(p)\times U(q))&&G(p,q)&&\Bbb P^{p-1}\times\Bbb P^{q-1}&&\text{Segre}&&&&\cr
\omit&height3pt&\omit&&\omit&&\omit&&\omit&&\omit&&\omit&&\cr
&&\text{II}&&SO(2n)&&U(n)&&G^{II}(n,n)&&G(2,n-2)&&\text{Pl\"ucker}&&&&\cr
\omit&height3pt&\omit&&\omit&&\omit&&\omit&&\omit&&\omit&&\cr
&&\text{III}&&Sp(n)&&U(n)&&G^{III}(n,n)&&\Bbb P^{n-1}&&\text{Veronese}&&&&\cr
\omit&height3pt&\omit&&\omit&&\omit&&\omit&&\omit&&\omit&&\cr
&&\text{IV}&&SO(n+2)&&SO(n)\times SO(2)&&Q^n&&Q^{n-2}&&\text{by $\Cal O(1)$}&&&&\cr
\omit&height3pt&\omit&&\omit&&\omit&&\omit&&\omit&&\omit&&\cr
&&\text{V}&&E_6&&\text{Spin}(10)\times U(1)&&\Bbb P^2(\Bbb O)\otimes_{\Bbb R}\Bbb C&&G^{II}(5,5)&&\text{by $\Cal O(1)$}&&&&\cr
\omit&height3pt&\omit&&\omit&&\omit&&\omit&&\omit&&\omit&&\cr
&&\text{VI}&&E_7&&E_6\times U(1)&&\text{exceptional}&&\Bbb P^2(\Bbb O)\otimes_{\Bbb R}\Bbb C&&\text{Severi}&&&&\cr
\omit&height3pt&\omit&&\omit&&\omit&&\omit&&\omit&&\omit&&\cr
\omit&\multispan{12}\hrulefill\cr
}}$$

\spb
It is possible to give an {\it $a$-priori\/} proof of the fact that
$\Cal C_o$ is Hermitian symmetric, as follows. Equipping $S=G/K$ with
a canonical K\"ahler-Einstein metric we can induce on
$\Bbb PT_o(S)$ and hence on $\Cal C_o$ a K\"ahler metric. In Mok [M2, App.(III.2),
Prop.1, pp.246ff] it was proven that holomorphic sectional curvatures
are pinched between $\frac 1{\phantom{,}2\phantom{,}}$ 
and 1 (up to normalization). By Ros [R],
we know that such K\"ahler submanifolds of the projective space
have to have parallel second fundamental forms and are hence in
particular Hermitian symmetric. It is remarkable that 
the list of standard cones exhausts the set of all such K\"ahler
submanifolds of projective spaces (equipped with the Fubini-Study metric).

\mpb\noindent
(2.2)\quad  Let now $\pi:X\to\triangle$, $X_t\cong S$ for $t\ne 0$,
be the regular family of compact K\"ahler manifolds as in Theorem 1.
Let $L\to X$ be the holomorphic line bundle such that $L|_{X_t}$ is a
positive generator of the Picard group $\text{Pic}(X_t)$ for all 
$t\in\triangle$. For $t\ne 0$, a degree-1 curve with respect to $L$
is necessarily rational and will be called a minimal rational curve.
Let $C\subset X_0$ be a limit of a sequence of minimal rational curves on
$\pi^{-1}(\triangle-\{0\})$ as cycles. Since $X_0$ is K\"ahler, and $C$ is of
degree 1 with respect to $L$, it must remain reduced and
irreducible (and rational). From now on a curve
$C\subset X_0$ will be called a minimal rational curve if and only if
it is a limit of minimal rational curves on $\pi^{-1}(\triangle-\{0\})$.

\spb
A holomorphic vector bundle $V$ over
a Riemann sphere $\Bbb P^1$ is said to be semipositive if and only if all
rank-1 factors in a Grothendieck decomposition of $V$ are
of degrees $\ge 0$. Let $f:\Bbb P^1\to X_0$ be a generically injective
holomorphic map, and write $C=f(\Bbb P^1)$.
By abuse of language we say that $C$ has semipositive normal bundle in $X_0$ if and only if
$f^*T_{X_0}$ is semipositive on $\Bbb P^1$, which is the same as saying
that $\text{Graph}(f)\subset \Bbb P^1\times X_0$ has semipositive normal
bundle as an embedded Riemann sphere.
Over $t\ne 0$ the Chow space of minimal rational curves fills up $X_t$.
It follows that there exists a unique irreducible component $\Cal K$ of the
Chow space of minimal rational curves on $X_0$, such that curves $C$
corresponding to points on $\Cal K$ sweep through $X_0$. 
As a consequence, there exists a subvariety $B\subset X_0$ such that
for $x\in X_0-B$ and for any
$[C]\in\Cal K$ passing through $x$, the normal bundle of $C$ in $X_0$ is
semipositive. 

\mpb\noindent
(2.3)\quad  We are going to construct quotient spaces of complex manifolds
under the action of complex Lie groups. 
This construction is well-known in the algebraic setting (cf. Mori [Mr]).
We will simply translate it to the complex-analytic setting.
Let $G$ be a complex Lie group
and $Z$ be a complex manifold. A holomorphic $G$-action on
$Z$ is by definition a holomorphic map $F:G\times Z\to Z$ such that 
(i) for any $g\in G$, $g(z):=F(g,z)$ is a biholomorphism of $Z$;
and (ii) $(gh)(z)=g(h(z))$ for any $g,h\in G$. We are interested to endow
the set of orbits $Z/G$ with a complex structure. For the formulation,
we say that the $G$-action is free if and only if $g(z)\ne z$ whenever
$g\ne e$, the identity element of $G$. We say that the $G$-action
has closed orbits if and only if $F(G\times\{z\})$ is closed in $Z$ 
for all $z\in Z$. The following lemma is standard.
(cf. Holmann [Ho, Satz 21]).

\spb
\proclaim{Lemma 1} Let $Z$ be a complex manifold, $G$ be a complex
Lie group and $F:G\times Z\to Z$ be a free holomorphic action
of $G$ on $Z$ with closed orbits. 
Then, the set of orbits $Z/G$ can be endowed the
structure of a complex manifold such that the canonical map 
$Z\to Z/G$ is a holomorphic submersion, realizing $Z$ as a holomorphic 
principal $G$-bundle.
\endproclaim

\spb
Consider the set $\Cal H$ of all holomorphic maps
$f:\Bbb P^1\to X$, such that $\deg(f^*L)=1$ and $f^*T$ is semipositive,
where $T\to X$ denotes the relative holomorphic tangent bundle
of $\pi:X\to\triangle$.
From the exact sequence $0\to f^*T \to f^*T_X\to\Cal O\to 0$ over
$\Bbb P^1$ and the vanishing of $H^1(\Bbb P^1,E)$ for semipositive
holomorphic vector bundles $E$, it follows that $f^*T_X\cong f^*T\oplus\Cal O$
whenever $f^*T\ge 0$, so that $f^*T_X\ge 0$. Conversely if
$f^*T_X\ge 0$, $f^*T_X\cong P\oplus\Cal O^\ell$ where $P$ is a
positive bundle. Then, $f^*T_X\to\Cal O$ is induced by a projection
$\Cal O^\ell\to\Cal O$, and the kernel $f^*T_X$ is isomorphic to
$P\oplus\Cal O^{\ell-1}$, hence semipositive. Thus, in the definition
of $\Cal H$, it is equivalent to assume $f^*T_X\ge 0$.
Furthermore, for $[f]\in\Cal H$, $f^*T_X\cong f^*T\oplus\Cal O$.

\spb
The complex Lie group $\text{\it Aut}(\Bbb P^1)\cong SL(2,\Bbb C)/\{\pm I\}$
acts on $\Cal H$ with closed orbits 
by the action $\varphi[f]=[f\circ\varphi^{-1}]$ for
any $[f]\in\Cal H$ and any $\varphi\in\text{\it Aut}(\Bbb P^1)$. Since
$f:\Bbb P^1\to X_0$ is generically injective the action is also free,
so that by Lemma 1 we have on the quotient space $\Cal D=\Cal H/\text{\it Aut}(\Bbb P^1)$
the structure of a complex manifold realizing $\Cal H\to\Cal H/\text{\it Aut}(\Bbb P^1)$
as a principal $\text{\it Aut}(\Bbb P^1)$-bundle. For each $[f]\in\Cal H$,
$f(\Bbb P^1)=C\subset X_t$ for some $t\in\triangle$, 
so that there is a canonical holomorphic map $\nu:\Cal H\to\triangle$.
The kernel of the differential $d\tau$ at $[f]\in\Cal H$ is given by
$\Gamma(\Bbb P^1,f^*T)\subset\Gamma(\Bbb P^1,f^*T_X)$, of corank 1,
so that $\nu:\Cal H\to\triangle$ is a holomorphic submersion.
Correspondingly we obtain a holomorphic submersion
$\Cal D\to\triangle$. The fiber $\Cal D_t$ over $t\ne 0$ can be
identified with the Chow space of minimal rational curves on $X_t$.
Over $t=0$, $\Cal D_0$ can be identified with the normalization of
a Zariski-open set of the irreducible component $\Cal K$ of the
Chow space of minimal rational curves on $X_0$.
Let now $\text{\it Aut}_0(\Bbb P^1)\subset\text{\it Aut}(\Bbb P^1)$
be the isotropy subgroup at $0\in\Bbb P^1$. Then,
$\Cal U:=\Cal H/\text{\it Aut}_0(\Bbb P^1)$ is a complex manifold.
As $\text{\it Aut}(\Bbb P^1)/\text{\it Aut}_0(\Bbb P^1)\cong\Bbb P^1$,
the canonical map $\tau: \Cal U\to\Cal D$ realizes $\Cal U$ as a
holomorphic $\Bbb P^1$-bundle over $\Cal D$. 

\spb
For $x\in X$ let $\Cal H_x\subset\Cal H$ be the subvariety consisting
of all $[f]\in\Cal H$ such that $f(0)=x$. We will say that $f$ is marked 
at $x$. We have

\spb
\proclaim{Proposition 4} For any $x\in X$, $\Cal H_x\subset\Cal H$ is a
complex submanifold. The canonical action of $\text{\it Aut}_0(\Bbb P^1)$
on $\Cal H_x$ induces on $\Cal M_x=\Cal H_x/\text{\it Aut}_0(\Bbb P^1)$
the structure of a complex manifold. Furthermore, the canonical projection
$\rho:\Cal U\mapsto X$ is a holomorphic submersion such that the fiber
over $x$ is $\Cal M_x$.
\endproclaim

\spb
\demo{Proof} For a holomorphic vector bundle $E$ on $\Bbb P^1$ we denote 
by $E(k)$ the twisted bundle $V\otimes\Cal O(k)$. Let $[f]\in\Cal H_x$.
Then, the deformation of $[f]$ within $\Cal H_x$ is unobstructed 
since $H^1(\Bbb P^1, f^*T_X(-1))=0$, $f^*T_X$ being semipositive.
Thus, $\Cal H_x\subset\Cal H$ is a (closed) complex submanifold such that
$T_{[f]}(\Cal H_x)=\Gamma(\Bbb P^1, f^*T_X(-1))\subset\Gamma(\Bbb P^1,
f^*T_X)=T_{[f]}(\Cal H)$. $\Gamma(\Bbb P^1, f^*T_X(-1))$ is of constant
$\text{rank}=\deg(f^* K_X^{-1})=\deg(K_X^{-1})$, where $K_X$
denotes the canonical line bundle over $X$.
Applying Lemma 1 to $\Cal H_x$ together with the canonical action of
$\text{\it Aut}_0(\Bbb P^1)$, we obtain on $\Cal M_x$ the structure of
a complex manifold. The canonical map $\Cal M_x=\Cal H_x/\text{\it Aut}_0(\Bbb P^1)
\to\Cal H/\text{\it Aut}_0(\Bbb P^1)=\Cal U$ identifies $\Cal M_x$ as a
complex submanifold. We denote by $\{f\}$ the class in $\Cal U$ defined by
$[f]$. Then $T_{\{f\}} (\Cal M_x)\cong\Gamma
(\Bbb P^1, f^*T_X(-1))/f_*\Gamma(\Bbb P^1, T_{\Bbb P^1}(-1))$.
Here $f_*$ denotes the homomorphism induced by the differential 
$df:T_{\Bbb P^1}\to f^*T_X$. The canonical map $\Cal H\to X$
sending $[f]$ to $f(0)$ induces a holomorphic map $\rho:\Cal U\to X$
such that $\text{Ker}\,d\rho(\{f\})=\Gamma(\Bbb P^1,f^*T_X(-1))/
f_*\Gamma(\Bbb P^1,T_{\Bbb P^1}(-1))=T_{\{f\}}(\Cal M_x)$, realizing
$\rho$ as a holomorphic submersion with fibers $\Cal M_x$, as desired.
\enddemo

\spb
By (2.2), there exists a subvariety $B\subset X_0$ such that for any
$[C]\in\Cal K$ passing through $x\in X_0-B$, the normal bundle of
$C$ in $X_0$ is semipositive. For such $x\in X_0-B$, by abuse of language
we will call the compact complex manifold $\Cal M_x$ the Chow space of
minimal rational curves marked at $x$. 
For the projection $\rho:\Cal U\to X$, we will sometimes write
$\Cal M$ for $\rho^{-1}(X-B)$. We have

\spb
\proclaim{Proposition 5} The canonical projection $\rho|_{\Cal M}:\Cal M\to
X-B$ is a regular family of projective-algebraic manifolds. Furthermore,
for $x\in X_t$, $t\ne 0$, the fiber $\Cal M_x$ is biholomorphic to the
standard cone $\Cal C_o$ of the model space $S$.
\endproclaim

\spb
\demo{Proof} The only thing to prove is the 
projective-algebraicity of $\Cal M_x$
for $x\in X_0-B$. Denote by $\Cal K_x\subset\Cal K$ the subvariety of
minimal rational curves passing through $x$. Since the canonical
map $\imath:\Cal D_0\to\Cal K$ is the normalization of a Zariski-open
subset of $\Cal K$ and the Chow space $\Cal K$ is projective-algebraic,
$\Cal D_0$ is quasi-projective. As $\tau|_{\Cal M_x}:\Cal M_x\to\Cal D_0$
is the normalization onto $\tau(\Cal M_x)$, $\Cal M_x$ is also
projective-algebraic, as desired.
\enddemo

\mpb\noindent
(2.4)\quad For the model space $S$ in Theorem 1 consider the first
canonical embedding $S\hookrightarrow\Bbb P^N$. For any minimal rational
curve $C\subset S$, which is a rational line in $\Bbb P^N$, $T_S|_C$
is a semipositive bundle which is a holomorphic subbundle of 
$T_{\Bbb P^N}|_C\cong\Cal O(2)\oplus [\Cal O(1)]^{N-1}$. 
It follows that $T_S|_C$ must be isomorphic to $\Cal O(2)\oplus[\Cal O(1)]^p
\oplus\Cal O^q$ for some nonnegative integers $p,q$. By considering
deformations of marked minimal rational curves, we conclude that $p$ is
the dimension of the standard cone $\Cal C_o$ of $S$.

\spb
Let $C\subset X_t$, $t\in\triangle$, be a minimal rational curve,
represented by $f:\Bbb P^1\to X$. We say that $C$ is standard if
and only if $f^*T_{X_t}\cong\Cal O(2)\oplus [\Cal O(1)]^p\oplus\Cal O^q$.
Thus $C$ is standard whenever $t\ne 0$. We note that on $X_0$, any
standard minimal rational curve must be immersed, since the only
non-trivial homomorphism in $\Gamma(\text{Hom}(T_{\Bbb P^1},f^*T_{X_0}))$
must be everywhere injective, as $T_{\Bbb P^1}\cong\Cal O(2)$. However,
we do not know {\it a priori\/} that standard minimal rational curves are
embedded in $X_0$. For the existence, we have

\spb
\proclaim{Proposition 6} On $X_0$ a generic minimal rational curve
with semipositive normal bundle is standard.
\endproclaim

\spb
For the proof of Proposition 6 we need the following lemma of Mori's.

\spb
\proclaim{Lemma 2 (Mori's Breaking-up Lemma ([Mr])} Let $\Sigma$ be an
algebraic surface and $\sigma:\Sigma\to\Gamma$ be a surjective holomorphic
map of $\Sigma$ over an algebraic curve $\Gamma$ such that the
generic fiber is irreducible and rational. Suppose there exist 
two distinct holomorphic sections $\Gamma_0$ and $\Gamma_\infty$
of $\Sigma$ over $\Gamma$, each of which can be blown down to a point.
Then, some fiber of $\sigma$ must be reducible.
\endproclaim

\spb
\demo{Proof of Proposition 6} The proof is similar to that of Mok 
([Mk2, Lemma (2.4.3), pp.203 ff.]). We argue by contradiction. 
If Proposition 6 fails, for every $[f]\in\Cal H$, $f(\Bbb P^1)$ is not standard.
At $[f]\in\Cal H$, $T_{[f]}(\Cal H)$ is identified with
$\Gamma(\Bbb P^1, f^*T_X)$. As all $f(\Bbb P^1)$ are non-standard,
we can choose a non-empty open subset $\Omega$ of $\Cal H$ and a
holomorphic vector field $\eta$ on $\Omega$ such that $\eta([f])\in
\Gamma(\Bbb P^1, f^*T_X\otimes [0]^{-1}\otimes [\infty]^{-1})
- f_*\Gamma(\Bbb P^1,T_{\Bbb P^1})$ for all 
$[f]\in\Omega$. Here for $\zeta\in\Bbb P^1$, $[\zeta]$ denotes
the positive divisor line bundle defined by $\zeta$. Fix
some $[f_0]\in\Omega$ such that $f_0(0)\ne f_0(\infty)$. Write
$x=f_0(0)$ and $y = f_0(\infty)$. An integral curve of $\eta$
in $\Cal H$ containing $[f_0]$ gives a non-trivial one parameter
family of minimal rational curves passing through $x$ and $y$.
As the Chow space $\Cal K$ is algebraic, this means that there exists
an algebraic family of minimal rational curves $C_\lambda$,
parametrized by some algebraic curve $\Gamma$ such that all $C_\lambda$
passes through the distinct points $x$ and $y$. Since all $C_\lambda$ are
irreducible this leads to a contradiction to Lemma 2. The proof of
Proposition 6 is completed.
\enddemo

\newpage

\bpb\noindent
{\bf \S 3 Non-deformability of normalized Chow spaces in the case of Grassmannians}

\mpb\noindent
(3.1) \quad Making use of the observation that standard cones are
themselves Hermitian symmetric spaces of the compact type, we formulate
the following proposition furnishing the inductive step in our proof
of Theorem 1.

\spb
\proclaim{Proposition 7} Let $S'$ be an irreducible  Hermitian symmetric 
manifold of the compact type of $\text{rank}\ge 2$ and of complex dimension
$n$. Suppose Theorem 1 has been established for all $S$ of complex
dimension $<n$. Then, for a regular family $\pi:X\to\triangle$ of compact
complex manifolds such that $X_t\cong S'$ for $t\ne 0$ and for a generic
point $x\in X_0$, the normalized 
Chow space $\Cal M_x$ of marked minimal rational
curves at $x$ is biholomorphic to the standard cone $\Cal C'_o$ of $S'$.
\endproclaim

\spb
For $S'$ different from Grassmannians, Proposition 7 follows from 
[(2.3), Proposition 4] on the smoothness of $\Cal M_x$ for generic $x\in X_0$.
For Grassmannians $S'\cong G(p,q)$, the standard cone $\Cal C'_o$ is
biholomorphic to $\Bbb P^{p-1}\times\Bbb P^{q-1}$, and we have to
worry about the possibility of a jump of complex structures in
deforming the latter manifold. Nonetheless the deformation we encounter is
subordinate to a deformation of Grassmannians, and we assert the following
proposition which implies Proposition 7.

\spb
\proclaim{Proposition 8} Let $\pi:X\to\triangle$ be a regular family
of compact K\"ahler manifolds and denote by $X_t$ the fiber
$\pi^{-1}(t)$ for all $t\in\triangle$. Suppose that for $t\ne 0$,
$X_t$ is biholomorphic to the Grassmannian $G(p,q)$, $p,q\ge 2$.
Let $s:\triangle\to X$ be a holomorphic section of $\pi:X\to\triangle$
such that the normalized Chow spaces of marked minimal rational curves 
$\{\Cal M_{s(t)}:t\in\triangle\}$ define a regular family 
$\tau:\Cal M\to\triangle$ of projective-algebraic 
manifolds, so that $\Cal M_{s(t)}$ is 
biholomorphic to $\Bbb P^{p-1}\times\Bbb P^{q-1}$ for $t\ne 0$.
Then, the central fiber
$\Cal M_{s(0)}$ is also biholomorphic to $\Bbb P^{p-1}\times\Bbb P^{q-1}$.
\endproclaim

\spb
In the example of the deformation of $\Bbb P^1\times\Bbb P^1$ to a
Hirzebruch surface $\Sigma_a$ with $a>0$ even (cf. Siu [S2, pp.97-98], 
for example),
limits of one of the two kinds of direct factors are decomposable.
We contend that such a phenomenon cannot occur in the case of
Proposition 8. We have 

\spb
\proclaim{Lemma 3} Let $\{t_i\}$ be a sequence of points on
$\triangle-\{0\}$ converging to 0. Suppose $Z_i\subset \Cal M_{t_i}$ are
projective spaces corresponding either to some 
$\{[a]\}\times\Bbb P^\ell$ or some $\Bbb P^k\times\{[b]\}$
under a biholomorphism
$\Cal M_{t_i}\cong\Bbb P^{p-1}\times\Bbb P^{q-1}$, where $1\le k\le p-1$
and $1\le\ell\le q-1$, and the projective spaces $\Bbb P^k$ resp.
$\Bbb P^\ell$ are understood to be projective-linear in $\Bbb P^{p-1}$
resp. $\Bbb P^{q-1}$. Suppose the sequence $\{Z_i\}$ converges as
cycles of $\Cal M$ to some $Z_0\subset\Cal M_0$. Then, $Z_0$ is
reduced and irreducible.
\endproclaim

\spb
\demo{Proof} We work first of all on
the model $S=G(p,q)$. Identifying $T_o(S)$, $o\in S$, with the tensor
product $\Bbb C^p\otimes\Bbb C^q$ in the standard way the cone
$\Cal C_o\subset T_o(S)$ is equivalently $\{[a\otimes b]:a\in\Bbb C^p$,
$b\in\Bbb C^q; a,b\ne 0\}$. Given $a\in\Bbb C^p$, the subvariety
$\{[a]\}\times\Bbb P^{q-1}\hookrightarrow\Bbb P^{p-1}\times\Bbb P^{q-1}
\hookrightarrow\Bbb P(T_0(S))$ is the projective-linear subspace
$\Bbb P(a\otimes\Bbb C^q)$. With respect to the Pl\"ucker
embedding $\nu:S\hookrightarrow\Bbb P^N$ minimal rational curves are
mapped biholomorphically onto rational lines in $\Bbb P^N$.
The projective-linear subspace $\Bbb P(a\otimes\Bbb C^q)\subset\Cal C_o$
then spans a $(q-1)$-dimensional family of rational lines in
$\Bbb P^N$ sweeping out a $q$-dimensional projective linear subspace
$\cong\Bbb P^q\hookrightarrow\Bbb P^N$. The same description is valid
for $\Bbb P(\Bbb C^p\otimes b)$.
\enddemo

\spb
Let now $\omega$ be the K\"ahler form of a Fubini-Study metric on
$\Bbb P^N$ such that the K\"ahler class $[\omega]$ 
is the positive generator of $H^2(\Bbb P^N,\Bbb Z)\cong\Bbb Z$.
For any $m$-dimensional projective linear subspace 
$\Bbb P^m\hookrightarrow\Bbb P^N$ we have $\int_{\Bbb P^m}\omega^m=1$.
The Pl\"ucker embedding $\nu:G(p,q)\to\Bbb P^N$ induces an
isomorphism $\nu^*:H^2(\Bbb P^N,\Bbb Z)\to H^2(G(p,q),\Bbb Z)\cong\Bbb Z$.
For $\{[a]\}\times\Bbb P^\ell$ resp. $\Bbb P^k\times\{[b]\}$ as
in the statement of Lemma 3, let $W\subset G(p,q)$ be the 
$(\ell+1)$ resp. $(k+1)$-dimensional subvariety swept out by the
corresponding set of marked minimal rational curves at 0.
Write $m$ for the dimension of $W$ (i.e. $k+1$ or $\ell+1$).
As the image $\nu(W)\subset\Bbb P^N$ is a projective linear
subspace we have
$$
\int_W (\nu^*\omega)^m = \int_{\nu(W)} \omega^m = 1\ .
$$
Let now $t_i\in\triangle$ and
$Z_i\subset\Cal M_{t_i}$ be as in the hypothesis of
Lemma 3 and denote by $W_i\subset X_{t_i}\cong G(p,q)$ the corresponding
$W$'s. The convergence of $Z_i$ to some $Z_0\subset\Cal M_0$ as a
subvariety implies the convergence of $W_i$ to some 
subvariety $W_0\subset X_0$. Since $W_i$ is of degree 1 with respect to
the ample line bundle $L$, $W_0$ is reduced and irreducible.
For $m=2$, i.e. $k$ resp. $\ell=1$, it follows that $Z_0$ (of
pure dimension 1) is reduced and irreducible, since any
irreducible component of $Z_0$ necessarily sweeps out a surface.
The general case can be deduced from the case of $k$ resp. $\ell=1$,
as follows. Suppose $Z_0$ is reducible. Pick $p,q\in Z_0$ to be
smooth points lying on different irreducible components, so that
any subvariety of $Z_0$ containing both $p$ and $q$ must be reducible.
Since $Z_i\subset\Cal M_{t_i}$ converges to $Z_0$ as cycles there 
exist pairwise distinct points $p_i,q_i\in Z_i$ such that $p_i\to p$
and $q_i\to q$ on $\Cal M$. Let $C_i\subset Z_i$ be the unique
projective line on $Z_i\cong\Bbb P^k$ resp. $\Bbb P^\ell$ containing
$p_i$ and $q_i$. Then, $C_i$ are of degree 1 with respect to $L$,
so that $C_i$ converges as cycles to a curve $C_o$ on $\Cal M_0$.
By the preceding argument for $k$ resp. $\ell=1$, $C_o$ must be
irreducible, contradicting with our choice of $p,q\in Z_0$.
We have thus proven that $Z_0$ is irreducible. It must then also be
reduced since $W_0$ is reduced.

\spb\noindent
(3.2) \quad It is possible to deduce Proposition 8 from Lemma 3
by an argument applicable to the more general situation of 
deformations of products of projective-algebraic manifolds.
We have chosen instead to present a proof adapted to the special
case of projective spaces more in line with methods used in the article,
using minimal rational curves and the tangent map.
Generic choices of $Z_0$ as in Lemma 3 will be shown to be
immersed projective spaces and the two families of such
cycles will be shown to intersect transversally. From this we will
deduce Proposition 8. 

\spb
Pick $t_0\in\triangle$, $t_0\ne 0$, and consider a cycle $Z'_{t_0}\subset
\Cal M_{t_0}$ corresponding to some $\Bbb P^{p-1}\times\{[b]\}$ under a
biholomorphism $\Cal M_{t_0}\cong\Bbb P^{p-1}\times\Bbb P^{q-1}$. Let
$\Cal D'$ be the irreducible component of the Chow space $\Cal D$
containing the point $[Z'_{t_0}]$ defined by $Z'_{t_0}$. 
Similarly choose $Z^{\prime\prime}_{t_0}\subset\Cal M_{t_0}$ corresponding
to some $\{[a]\}\times\Bbb P^{q-1}$ and define $\Cal D''\subset
\Cal D$ to be the corresponding irreducible component of $\Cal D$.
We have canonical projection maps $\rho':\Cal D'\to\triangle$ and
$\rho'':\Cal D''\to\triangle$. Write $\Cal D'_t$ resp. $\Cal D''_t$ for the
fibers over $t\in\triangle$. 
By studying rational curves on $Z'$, $[Z']\in\Cal D'$;
resp. $Z''$, $[Z'']\in\Cal D''$; we prove

\spb
\proclaim{Lemma 4} There exists a non-empty Zariski-open subset
$U=\Cal M_0-E\subset \Cal M_0$ such that for every 
$[Z'_0]\in\Cal D'_0$, $Z'_0$ is an immersed projective space of
complex dimension $p-1$ provided that $Z'_0\cap U\ne\emptyset$. 
Furthermore, for every 
$[Z''_0]\in\Cal D''_0$, $Z''_0$ is an immersed projective space of
complex dimension $q-1$ provided that $Z''_0\cap U\ne\emptyset$.
\endproclaim

\spb
\demo{Proof} First of all, we observe that the biholomorphisms
$\Cal M_t\cong\Bbb P^{p-1}\times\Bbb P^{q-1}$ for $t\ne 0$ define
integrable holomorphic subbundles $T'_t$, $T''_t\subset T_t:=T_{\Cal M_t}$
on $\Cal M_t$ such that the direct factors $Z'_t\subset \Cal M_t$
corresponding to $\Bbb P^{p-1}\times\{[b]\}$, resp.
$Z''_t\subset\Cal M_t$ corresponding to $\{[a]\}\times\Bbb P^{q-1}$,
are integral submanifolds of $T'_t$ resp. $T''_t$. 
As the distribution $z\mapsto T'_t(z)$ resp. $z\mapsto T''_t(z)$ on
$\Cal M_t$ are defined by holomorphic maps into Grassmannians,
they extend meromorphically
to $\Cal M$, so that there exists a proper complex-analytic subvariety
$E_1\subset\Cal M_0$ and holomorphic vector subbundles $T'$,
$T''\subset T$ on $\Cal M-E_1$ such that for $t\ne 0$, $T'\big|_{\Cal M_t}=T'_t$,
$T''\big|_{\Cal M_t}=T''_t$. We now assume that $[Z'_0]\in\Cal D'_0$ and
$y\in Z'_0$ have been so chosen that $y\notin E_1$.
\enddemo

\spb
Pick $t_0\ne 0$ and $Z'_{t_0}\subset\Cal M_{t_0}$ such that
$[Z'_{t_0}]\in\Cal D'_{t_0}$. We have $Z'_{t_0}\cong\Bbb P^{p-1}$.
Let $C_{t_0}$ be a minimal rational curve on $Z'_{t_0}\cong\Bbb P^{p-1}$.
For the Chow space $\Cal D$ of $\Cal M$, let $\Cal E'\subset\Cal D$ be the
irreducible component containing the point $[C_{t_0}]$.
From now on, we replace $\Cal E'$ by its normalization and
denote by $\Cal E'_t$ the fiber over $t$ of the
canonical projection $\Cal E'\mapsto\triangle$. Similar notations
will be used for canonically fibered spaces over $\triangle$. 
From Lemma 3 we know that for every $[C_0]\in\Cal E'_0$, 
$C_0$ is an irreducible rational curve, represented by 
$f_0:\Bbb P^1\to\Cal M_0$. There exists a unique irreducible
component $\Cal E'_{01}$ of $\Cal E'_0$ whose members cover $\Cal M_0$.
Members of all other irreducible components cover a proper subvariety 
$E_2\subset\Cal M_0$. By Lemma 3 and the Breaking-up Lemma of Mori,
we know that for a generic choice of $C_0\subset\Cal M_0$
with $[C_0]\in\Cal E'_{01}$, represented by $f_0:\Bbb P^1\to\Cal M_0$,
we have
$$
f_0^* T_{\Cal M_0} \cong \Cal O(2) \oplus [\Cal O(1)]^{p-2} \oplus\Cal O^{q-1}\ .
\tag $*$
$$
In analogy with (2.4) $C$ is said to be standard whenever $t\ne 0$; or
$t=0$ and $(*)$ is valid for 
$f_0:\Bbb P^1\to\Cal M_0$ defining $C$.  We have

\spb\noindent
{\smc Claim:}\quad For a generic choice of $w_0$ on $\Cal M_0$,
every rational curve $C_0\subset\Cal M_0$ with $[C_0]\in\Cal E'_0$,
$w_0\in C_0$, is standard.

\spb
\demo{Proof of Claim} 
Since $\cup\{C_0:[C_0]\in\Cal E'_{01}\}$ covers $\Cal M_0$
there exists a proper complex-analytic subvariety $E_3\subset\Cal M_0$
such that for any $w_0\in\Cal M_0-E_3$ and any $[C_0]\in\Cal E'_{01}$
passing through $w_0$, defined by $f:\Bbb P^1\to\Cal M_0$, we have
$f^*T_{\Cal M_0}\ge 0$. 
The semipositivity of $f^*T_{\Cal M_0}$ implies the smoothness of
$\Cal E'_{01}$ at $[C_0]$. Write $E=E_1\cup E_2\cup E_3$.
Suppose furthermore that $w_0\notin E$ and $C_0$ is smooth at $w_0$.
We are going to prove that $C_0$ is standard. 
Denote by $\Cal U\subset\Cal E'\times\Cal M$ the subvariety 
$\{([C],w)\in\Cal E'\times\Cal M: w\in C\}$. For any $t\in\triangle$,
$[C_t]\in\Cal E'_t$ and any smooth point $w\in C_t$ we define
$\Theta([C_t],w)\in\Bbb PT_w(\Cal M_t)$ to be the projectivization
of $T_w(C_t)$. We call $\Theta$ the tangent map. Pick a smooth
neighborhood $V$ of $([C_0],w_0)$ in $\Cal U$ with the following
properties: (i) $([C],w)\in V$ implies that $w\notin E$ and $w$ is
smooth on $C$; (ii) the canonical projection $\mu:V\to\Cal M$ has
constant rank and has connected fibers. We note that 
$\Theta([C],w)\in\Bbb PT'_w$ for $([C],w)\in V$. This is the case 
by the definition of $T'$ whenever $C$ lies over $t\ne 0$.
It remains true for $C$ lying over 0 by continuity, as $T'\subset T$
is a holomorphic vector subbundle on $\Cal M-E_1$. We note also that
both $V$ and $\Bbb PT'$ are of complex dimension $2p+q-3$. Since
$\dim_{\Bbb C}V = \dim_{\Bbb C}\Bbb PT'$, the tangent map
$\Theta:V\to\Bbb PT'$ is either an open immersion or it is
unramified exactly outside a non-trivial divisor $D\subset V$.
As any $[C_t]\in\Cal E'_t$
represents a standard minimal rational curve $C_t\subset\Cal M_t$ whenever
$t\ne 0$, the latter possibility can occur only if all $C_0$,
$[C_0]\in\Cal E'_0\cap V$ are not standard, contradicting the existence
of standard rational curves in $\Cal E'_{01}$.
We have proven that for $w_0\in\Cal M_0-E$, and for any $[C_0]\in
\Cal E'_{01}$ containing $w_0$ as a smooth point, $C_0$ is standard.
The same can be said for any $[C_0]\in\Cal E'_{01}$ passing through $w_0$. 
In fact, $C_0$ must contain a smooth point $z\in\Cal M_0-E$.
Repeating the same argument at $z$ in place of $w_0$ we conclude  that
$C_0$ must be standard.  The proof of the claim is completed.
\enddemo

\spb
Let now $w_0\in\Cal M_0-E$ and 
$[Z'_0]\in\Cal D'_0$ be such that $w_0\in Z'_0$. 
The holomorphic subbundle
$T'_0\subset T_0$ is defined at $w_0$. For any $[C_0]\in\Cal E'_0$
represented by $f_0:\Bbb P^1\hookrightarrow\Cal M_0$, $f_0(0)=w_0$, we know
that $f_0$ is an immersion at 0, so that $df_0(0)$ sends $T_0(\Bbb P^1)$
to a complex line $L_{C_0}$ in $T'_{w_0}$. Consider the normalized Chow
space $\Cal C_{w_0}^*$ of rational curves marked at $w_0$ and denote by
$\Cal C'_{w_0}$ the irreducible component of $\Cal C_{w_0}^*$ containing 
$[C_0]$. Then $\Cal C'_{w_0}$ is a $(p-2)$-dimensional complex manifold
such that the tangent map $\Theta:\Cal C'_{w_0}\to\Bbb PT'_{w_0}\cong\Bbb P^{p-2}$
is a local biholomorphism and thus a biholomorphism. The subvariety
$S=\cup\{C_0:[C_0]\in\Cal C'_{w_0}\}\subset\Cal M_0$ must have a
tangent space  = $\cup\{L_{C_0}:[C_0]\in\Cal C'_{w_0}\}=T'_{w_0}$.
For the proof of Lemma 4 we may assume that $p\ge 3$, as
the case $p=2$ is already implied by the Claim.
The family of rational curves parametrized by $\Cal C'_{w_0}$ gives
rise to a $\Bbb P^1$-bundle $B$ over $\Cal C'_{w_0}\cong\Bbb P^{p-2}$
with a distinguished section $\Sigma$ corresponding to the base point ${w_0}$.
The fact that $\Theta$ is a biholomorphism implies that the normal bundle
of $\Sigma\cong\Bbb P^{p-2}$ in $B$ is isomorphic to $\Cal O(-1)$. For $p\ge 3$
one can then blow down $\Sigma$ to a smooth point $b$ to obtain $\check B\cong\Bbb P^{p-1}$.
The canonical map $B\mapsto\Cal M_0$ induces a holomorphic map 
$\nu:\Bbb P^{p-1}\cong\check B\to\Cal M_0$. Since each
$[C_0]\in\Cal C'_{w_0}$ represents a standard rational curve, 
$\nu:\Bbb P^{p-1}\to\Cal M_0$ is an immersion outside
$b\in\check B$. Since $\nu(\Bbb P^{p-1})=S$ is smooth at ${w_0}$ and $\nu$ is
an immersion outside of a single point 0, of codimension $p-1\ge 2$,
$\nu$ must be a local biholomorphism at 0.
The same argument applies {\it verbatim\/} to $\Cal D''$ in place of 
$\Cal D'$, proving Lemma 4.

\spb
In the proof of Lemma 4 the integrable holomorphic subbundles 
$T'$, $T''\subset T$ were defined on $\Cal M-E$. For $x\in\Cal M-E$
denote by $Z'(x)$ (resp. $Z''(x)$) the unique leaf of the integrable
holomorphic distribution $x\mapsto T'_x$ (resp. $T''_x$) passing through
$x$. For $t\ne 0$ and $x\in\Cal M_t\cong\Bbb P^{p-1}\times\Bbb P^{q-1}$,
$Z'(x)\cap Z''(x)=\{x\}$, the intersection being transversal at $x$.
For the proof of Proposition 8 we need the following lemma extending
this property across generic points of $\Cal M_0$.

\spb
\proclaim{Lemma 5} In the notations of Lemma 4, for any $x\in\Cal M_0-E$,
$Z'(x)$ and $Z''(x)$ are smooth at $x$. Furthermore, 
$Z'(x)\cap Z''(x)=\{x\}$, the intersection being transversal at $x$.
\endproclaim

\spb
\demo{Proof} Consider a holomorphic family of embeddings
$\nu'_t:\Bbb P^{p-1}\to\Cal M_t$, $t\in\triangle$, with 
$\nu'_0=\nu'$; $\nu'(\Bbb P^{p-1})=Z'(x)$ and similarly a holomorphic
family of embeddings $\nu''_t:\Bbb P^{q-1}\to\Cal M_t$, $t\in\triangle$,
with $\nu''_0=\nu''$, $\nu''(\Bbb P^{q-1})=Z''(x)$. We know that for $t\ne 0$
$$
\nu'_t(\Bbb P^{p-1}) \cap \nu''_t(\Bbb P^{q-1}) = \{x_t\}
$$
for some $x_t\in\Cal M_t$, such that $x_t\mapsto x$. Write $[\Bbb P^{p-1}]$
resp. $[\Bbb P^{q-1}]$ for the fundamental class of $\Bbb P^{p-1}$ resp.
$\Bbb P^{q-1}$. We have
$$
(\nu'_t)_* [\Bbb P^{p-1}] \cdot (\nu''_t)_* [\Bbb P^{q-1}] = 1
$$
for $t\ne 0$ and hence also for $t=0$. Suppose $\nu'_0(\Bbb P^{p-1})\cap
\nu''_0(\Bbb P^{q-1})$ is zero-dimensional. Then, $Z'(x)$ and $Z''(x)$ must
be smooth at $x$, $x$ is the unique 
point of intersection and the intersection must be transversal, 
as desired. It remains therefore to rule out the existence
of a positive-dimensional irreducible component $S$ in $\nu'_0(\Bbb P^{p-1})
\cap\nu''_0(\Bbb P^{q-1})$. Given any such $S$ of complex dimension
$s$ with fundamental class $[S]$,
$$
0 \ne [S]\in(\nu'_0)_* H_{2s} (\Bbb P^{p-1},\Bbb Z) \cap
(\nu''_0)_* H_{2s} (\Bbb P^{q-1},\Bbb Z)\ . \tag 1
$$
Define $S'_0=(\nu'_0)^{-1}(S)\subset\Bbb P^{p-1}$, 
$S''_0=(\nu''_0)^{-1}(S)\subset\Bbb P^{q-1}$. Since
$\nu'_0:\Bbb P^{p-1}\to\Cal M_0$ and $\nu''_0:\Bbb P^{q-1}\to\Cal M_0$
are immersions, $S'_0$ and $S''_0$ are of pure dimension $s$.
Define $S'_t=\nu'_t(S'_0)$, $S''_t=\nu'_t(S''_0)$.
On the one hand,
$$
(\nu'_t)_* H_{2s} (\Bbb P^{p-1},\Bbb Q) \cap (\nu''_t)_* H_{2s}
(\Bbb P^{p-1},\Bbb Q) = \{0\} \ .\tag 2
$$
On the other hand, we have $[S'_t]\to k'[S]$; $[S''_t]\to k''[S]$
as homology classes for some positive integer $k'$ and $k''$. This
leads plainly to a contradiction between (1) and (2). The proof of
Lemma 5 is completed.
\enddemo

\spb
We are now ready to give a proof of Proposition 8.

\spb
\demo{Proof of Proposition 8} Recall that $\tau:\Cal M\to\triangle$
is a regular family such that $\Cal M_t=\tau^{-1}(t)$ is
biholomorphic to $\Bbb P^{p-1}\times\Bbb P^{q-1}$ for $t\ne 0$.
Denote by $T$ the relative holomorphic tangent bundle of 
$\tau:\Cal M\to\triangle$. By Lemma 5, there exists a proper complex-analytic
subvariety $E\subset\Cal M_0$, holomorphic vector subbundles
$T'$ and $T''$ of $T\big|_{\Cal M-E}$,
such that $T'_x\cap T''_x = \{0\}$, $T_x=T'_x\oplus T''_x$ for
any $x\in\Cal M-E$; and such that for $x\in\Cal M_t$, $t\ne 0$, the direct sum 
decomposition $T_x = T'_x\oplus T''_x$ arises from the product
decomposition $\Cal M_t \cong \Bbb P^{p-1}\times\Bbb P^{q-1}$. Given a
complex vector space $W$ of complex dimension $N$, positive integers
$N'$, $N''$ such that $N=N'+N''$, the set of direct sum decompositions
$W=W'\oplus W''$, with $W'$ resp. $W''\subset W$ of complex dimensions
$N'$ resp. $N''$, is parametrized by $\Cal S = GL(N,\Bbb C)/GL(N',\Bbb C)
\times GL(N'',\Bbb C)$. Since $G=GL(N',\Bbb C)\times GL(N'',\Bbb C)$ is
reductive and $E\subset\Cal M$ is a subvariety of
codimension $\ge 2$, we can apply the extension result of
[(1.3), Proposition 1] to extend the $G$-structure; i.e.,
equivalently the direct sum decompositions $W=W'\oplus W''$.
As a consequence, we obtain on all of $\Cal M$ a decomposition
$T=T'\oplus T''$ as a direct sum of holomorphic vector subbundles.
\enddemo

\spb
For $t\ne 0$, the holomorphic distributions on $\Cal M_t$ defined by
$x\mapsto T'_x$, $T''_x$ are integrable with closed leaves. The same is
true at $t=0$ by taking limits as $t\mapsto 0$. On $\Cal M_0$, the
holomorphic distribution $x\mapsto T'_x$ (resp. $x\mapsto T''_x$)
defines a regular family of compact complex manifolds such that the
generic fiber is biholomorphic to $\Bbb P^{p-1}$ (resp. $\Bbb P^{q-1}$),
by Lemma 4. $Z'(x)\cap Z''(x)=\{x\}$ remains valid by
considering intersection numbers (as in the proof of Lemma 4).
From this it follows that $\Cal M_0 \cong \Bbb P^{p-1}\times\Bbb P^{q-1}$,
as desired.

\spb\noindent
{\smc Remarks}

\noindent
For the Steinness of $\Cal S = GL(N,\Bbb C)/GL(N',\Bbb C)
\times GL(N'',\Bbb C)$, in place of applying Matsushima-Morimoto [MM],
we can use the explicit
description of $\Cal S$, as the complement of a {\it hypersurface\/}
$D$ in a product $G'\times G''$ of two Grassmannians, $G'=G(N',N'')$,
$G''=G(N'',N')$, where $D$ is of bidegree $(a,b)$, $a,b>0$, arising
from taking determinants of a set of $N$ vectors of $W\cong\Bbb C^N$.
$\Cal S = (G'\times G'')-D$ implies that $\Cal S$ is in fact affine-algebraic.
In particular, $\Cal S$ is a Stein manifold.

\bpb\noindent
{\bf \S 4 The linear span of tangents to minimal rational curves in Fano}
\linebreak
\phantom{\bf \S 4} {\bf manifolds}

\mpb\noindent
(4.1)\quad For the proof of Theorems 1 and 1$'$ we are going to consider
distributions arising from cones of minimal rational curves
in the central fiber $X_0$. The key question is integrability. In this
section we will formulate results on integrability and non-integrability
in the broader context of Fano manifolds.

\spb
Let $M$ be a Fano manifold. By Mori [Mr], $M$ is uniruled. Let $E$ be
the set of all irreducible components $\Cal P$ of the Chow space of
$M$ such that generic points of $\Cal P$ correspond to rational
curves with semipositive normal bundles. Fix a positive line bundle $L$ on
$M$ and define $\delta_L(\Cal P)$ as the degree of members of $\Cal P$
with respect to $L$. From the proof of  [(2.4), Proposition 6],
minimizing $\delta_L(\Cal P)$ among $\Cal P\in E$ we obtain 
$\Cal K\in E$ with the following property: A generic point of
$\Cal K$ corresponds to an immersed rational curve $C$ given by
$f:\Bbb P^1\to M$ such that $(*) f^*T_{M}\cong\Cal O(2)\oplus
[\Cal O(1)]^p\oplus\Cal O^q$ for some $p,q\ge 0$. (We note that
the arguments in [(2.4), Proposition 6] apply because we can choose
points $x,y\in M$ such that every rational curve passing through
$x$ resp. $y$ has semipositive normal bundle.) There may be several
possible choices of $\Cal K$. From now on we fix one such $\Cal K$.
In analogy with (2.4) we will call $C$ a standard minimal rational
curve if and only if (implicitly) $[C]\in\Cal K$ and the splitting condition
$(*)$ is satisfied.

\spb
Let $\Cal R\subset\Cal K$ be the Zarishi-dense subset corresponding
to standard minimal rational curves. Consider the universal family
maps $\psi:\Cal U\to\Cal R$, $\varphi:\Cal U\to M$. For a given
point $x\in M$, let $\Cal R_x=\psi(\varphi^{-1}(x))$ be the
set of standard minimal rational curves containing $x$. From now on
to streamline the discussion we will make the following assumption,
which is valid in the situation of the present article:

\spb
(Assumption) For a generic point $x\in M$, $\Cal R_x$ is irreducible.

\spb
The tangent map $\Theta_x:\Cal R_x\to\Bbb P T_x(M)$, which sends
$r\in\Cal R_x$ to the tangent vector of $\varphi(\psi^{-1}(r))$ at $x$,
is a holomorphic map. Let $\Cal K_x$ be the closure of $\Cal R_x$.
Then $\Theta_x$ induces a dominant rational map $\Cal K_x\to\Cal C_x$,
where $\Cal C_x$ is the closure of the image $\Theta_x(\Cal R_x)$.
Consider the holomorphic distribution defined on a Zariski-dense
subset of $M$ by the linear span of $\Cal C_x$ in $\Bbb P T_x(M)$.
We can extend this distribution to a distribution $W$ on $M$ outside
a subvariety $S(W)$ of codimension $>1$. Of course, $W$ is the
trivial distribution of the whole tangent bundle, if $\Cal C_x$ is
linearly non-degenerate in $\Bbb PT_x(M)$. We will sometimes refer to
$W$ as a meromorphic distribution on $M$ and say that $W$ is
integrable to mean that it is integrable on $M-S(W)$. 

\spb
At a generic point $x\in M$ let $W_x\subset T_x(M)$ be the subspace
defining $W$. By definition, $\Cal C_x\subset\Bbb PW_x$ is linearly
non-degenerate. A line in $\Bbb PW_x$ corresponds to a 2-plane in $W$,
defining a point in $\Bbb P\Lambda^2 W_x$. Consider the tangential 
lines to $\Cal C_x$ at smooth points of $\Cal C_x$. The closure of the
corresponding points on $\Bbb P\Lambda^2 W_x$ defines the variety of
tangential lines $\Cal T_x\subset\Bbb P\Lambda^2 W_x$.

\spb\noindent
(4.2)\quad We develop now a sufficient condition for the integrability
of the meromorphic distribution $W$. We have 

\spb
\proclaim{Proposition 9} Suppose the variety of tangent lines 
$\Cal T_x\subset\Bbb P\Lambda^2 W_x$ is linearly non-degenerate at
a generic point $x\in M$. Then, $W$ is integrable.
\endproclaim

\spb
As $W$ is defined algebro-geometrically using minimal rational curves we
do not have a straight-forward way to verify the Frobenius condition
on $W$. Instead, we will first of all reduce the verification to finding
at generic points $x\in X_0$ a set of integral surfaces $\Sigma$
such that $\Lambda^2 T_x(\Sigma)$ generates $\Lambda^2 W_x$.
Then, we will show that such integral surfaces $\Sigma$ can be 
constructed as pencils of rational curves.

\spb
We start with the following general statement on the Frobenius condition.

\proclaim{Lemma 6} Let $\Omega\subset\Bbb C^n$ be a domain and 
$D\subset T_\Omega$ be a holomorphic distribution on $\Omega$. 
Then, $D$ is integrable if and only if the following holds.

\spb
\item{$(*)$} Given $x\in\Omega$, there exist $D$-valued 
holomorphic vector fields
$\alpha_j$, $\beta_j$ defined on a neighborhood of $x$,
$1\le j\le N$ for some positive integer $N$, such that

\spb
\itemitem{\rm (i)} $[\alpha_j,\beta_j](x)\in D_x$,

\spb
\itemitem{\rm(ii)} $\{\alpha_j(x)\land\beta_j(x)\}$ spans $\Lambda^2 D_x$.
\endproclaim

\spb
\demo{Proof} The ``only if'' part is obvious. For the ``if'' part by
Frobenius we need to prove that $[D,D]\subset D$.
Fix $x\in\Omega$. Let $A$, $B$ be holomorphic vector fields with values in $D$
on a neighborhood $U$ of $x$, and $f$, $g$ be holomorphic functions on $U$. 
Then,
$$
[fA,gB] = fg[A,B] + fA(g)B - gB(f)A \equiv fg[A,B]\ \text{mod}\,D\ .
$$
It follows that the Lie bracket defines at $x$ a skew-symmetric
bilinear form on $D_x$ with values in the quotient space $T_x/D_x$.
The hypotheses (i) and (ii) then imply that this bilinear form vanishes
for every $x\in\Omega$. In other words $[D,D]\subset D$ and $D$ is
integrable, as desired.
\enddemo

\spb
For the proof of Proposition 9 the crucial point is to verify the
hypothesis of Lemma 6 for $W$ on some Zariski-dense subset
of $M$ using the deformation theory of rational curves.
For this purpose we have

\spb
\proclaim{Proposition 10} Let $x\in M$ be a generic point and
$C\subset M$ be a standard minimal rational curve passing 
through $x$, represented by $f:\Bbb P^1\to M$ such that $f(0)=x$. Write
$$
f^*T\big|_{\Bbb P^1} \cong\Cal O(2) \oplus [\Cal O(1)]^p \oplus\Cal O^q
$$
for a decomposition of $f^*T$ into a direct sum of line bundles. 
Write $df(T_0(\Bbb P^1))=\Bbb C\alpha$ and denote by
$P_\alpha\subset T_x$ the vector subspace corresponding to
$(\Cal O(2)\oplus [\Cal O(1)]^p)_0$. Let $\xi\in P_\alpha$ be a vector
such that $\alpha$ and $\xi$ are linearly independent. Then, there exists
a smooth complex-analytic surface $\Sigma$ on some neighborhood $\Omega$
of $x$ such that 
$$
T_x(\Sigma) \cong\Bbb C\alpha + \Bbb C\xi
$$
and such that at every point $y\in\Omega$, $T_y(\Sigma)\subset W_y$.
\endproclaim

\spb
\demo{Proof} To simplify notations
we assume that $C$ is embedded. Choose a point $y\in C\cap\Omega$,
$y\ne x$. Let $s\in\Gamma(C,T_X)$ be a holomorphic section such that
$s(x)=\xi$ and $s(y)=0$. This is possible because $\xi\in P_\alpha$.
Let now $\{C_\zeta: \zeta\in\triangle\}$ be a family of smooth standard
minimal rational curves given by $f_\zeta:\Bbb P^1\to X_0$ such that
$f_0\equiv f$ and such that $f_\zeta(\infty)=y$ for all $\zeta\in\triangle$.
Then, for any $\zeta_0\in\triangle$ consider the section
$s=\frac{\partial f_\zeta}{\partial\zeta}\Big|_{\zeta=\zeta_0}\in
\Gamma(\Bbb P^1,f^*_{\zeta_0} T_X\otimes\Cal O(-1))$. We have
$s(\infty)=0$. Identify $s$ as a section in $\Gamma(C_{\zeta_0};T_X)$.
For any $u\in C_{\zeta_0}$ with $T_u(C_{\zeta_0})=\Bbb C\alpha'$ we must
have $s(u)\in P_{\alpha'}$ as $s(y)=0$. Let now $\Sigma$ be the
image of $\{(z;\zeta): |z|<1; |\zeta|<\varepsilon\}$
under the map $F:\Bbb P^1\times\triangle\to X_0$ defined by 
$F(z,\zeta)=f_\zeta(z)$. We may assume that $\Sigma\subset\Omega$.
Since $\alpha\land\xi\ne 0$, 
for $\varepsilon$ sufficiently small $\Sigma$ is a locally closed
complex-analytic surface containing $x$ such that for any point
$u\in\Sigma$, $u=F(z,\zeta_0)$,
$$
T_u(\Sigma) = \Bbb C\alpha' + \Bbb Cs(u) \subset P_{\alpha'}\subset W_u
$$
for $s$ corresponding to 
$\frac{\partial f}{\partial\zeta}\Big|_{\zeta=\zeta_0}$.
The proof of Proposition 10 is completed.
\enddemo

\spb
We proceed with the proof of Proposition 9.

\spb
\demo{Proof of Proposition 9} At a generic point of $M$, $W_x$ 
is defined and is spanned by vectors $\alpha$ tangent to minimal
rational curves. By Proposition 10, given such an $\alpha$ and
$\xi\in P_\alpha$, $\alpha\wedge\xi\ne 0$, there exists an integral
surface $\Sigma$ of $W$ in a neighborhood of $x$, such that
$T_x(\Sigma)=\Bbb C\alpha\oplus\Bbb C\xi$. Thus there exist
$W$-valued holomorphic vector fields $\tilde\alpha$, $\tilde\xi$
on a neighborhood of $x$ such that $\tilde\alpha(x)=\alpha$;
$\tilde\xi(x)=\xi$ and $[\tilde\alpha,\tilde\xi](x)=0$. By assumption
the linear span of such $\alpha\wedge\xi$ at $x$ is 
$\Lambda^2W$, so that the hypothesis of Lemma 6 is satisifed on
a Zariski-dense subset of $M-\Cal S(W)$. It follows that $W$ is integrable,
as desired.
\enddemo

\mpb\noindent
(4.3)\quad We proceed to examine consequences of integrability
of the meromorphic distribution $W$ defined by rational curves.
When cones $\Cal C_x\subset\Bbb PT_x$ are linearly degenerate,
i.e., when $W$ is non-trivial, we will show that $b_2(M)\ge 2$.
To start with, we have

\spb
\proclaim{Proposition 11} 
Suppose the meromorphic distribution $W$ is integrable. Then every
leaf ${\Cal W}$ on $M-{\Cal S}(W)$ is closed and its topological closure
$\overline{\Cal W}$ is a complex-analytic subvariety.
\endproclaim

\spb
\demo{Proof} Let $x \in M$ be a generic point. It is enough 
to show that the leaf through $x$ can be compactified to a subvariety.

\spb
Let $\psi: {\Cal U} \rightarrow {\Cal R}$  and 
$\varphi: {\Cal U} \rightarrow M$
be the universal family. We know that there exists a subvariety
$Z \subset M$ so that ${\Cal R}_{x}$ and ${\Cal C}_{x}$ are
irreducible for $x \in M- Z$.

\spb
Consider an irreducible subvariety $A \subset M$ which is not contained in
${\Cal S}(W) \cup Z$ and whose generic part is contained in 
a single leaf of $W$. We will define a new irreducible variety ${\Cal V}(A)$ 
with the properties: (i) ${\Cal V}(A)$
is not contained in ${\Cal S}(W) \cup Z$
and (ii) the  generic part of ${\Cal V}(A)$ is contained in 
the leaf of $W$ containing the generic part of $A$.
 
\spb
The variety  $\psi(\varphi^{-1}(A))$ may
be reducible. But there is an irreducible component ${\Cal R}(A)$
which contains the irreducible ${\Cal R}_{a}$ 
for a generic $a \in A$. Let ${\Cal V}(A)$ be the closure of $\varphi
(\psi^{-1}({\Cal R}(A)))$. Then ${\Cal V}(A)$ is also irreducible and it 
contains $A$. From this, Property (i) is immediate.
To check Property (ii), note that
generic points of any curve $C$ with $[C] \in {\Cal R}$
which is not contained in ${\Cal S}(W)$, 
is contained in a single leaf of $W$. Curves corresponding to generic points
of ${\Cal R}(A)$ intersect $A-{\Cal S}(W)$. Hence generic points
of  such curves are
contained in the leaf of $W$ containing $A$. But such curves are dense in
${\Cal V}(A)$. It follows that generic points of ${\Cal V}(A)$ are contained
in the leaf of $W$ containing the generic part of $A$.

\spb
Now given a generic point $x \in M$, consider the irreducible varieties
${\Cal V}^{i}(x)$ defined inductively by ${\Cal V}^{1}(x) = {\Cal V}(x)$ and
${\Cal V}^{i+1}(x) = {\Cal V}({\Cal V}^{i}(x))$. 
Since ${\Cal V}^{i}(x) \subset {\Cal V}^{i+1}(x)$,
and they are irreducible, we get ${\Cal V}^{n}(x) = {\Cal V}^{n+1}(x)$, where
$n$ is the dimension of $M$. This implies that for a generic smooth
point $y \in {\Cal V}^{n}(x)$, a  curve corresponding to a generic point 
of ${\Cal R}_{y}$ is contained in ${\Cal V}^{n}(x)$. Hence ${\Cal C}_{y}
\subset {\Bbb P}T_{y}({\Cal V}^{n}(x))$, 
which implies that  $W_{y} \subset T_{y}({\Cal V}^{n}(x))$. 
But the generic part of ${\Cal V}^{n}(x)$ is  
contained in the leaf of $W$ containing $x$,
which implies $T_{y}({\Cal V}^{n}(x)) \subset W_{y}$. It follows that at a
generic point, ${\Cal V}^{n}(x)$ is the leaf of the foliation $W$. 
Thus the leaf of  $W$ containing $x$
can be compactified to the subvariety ${\Cal V}^{n}(x)$.
\enddemo

\spb
\proclaim{Proposition 12} Suppose the meromorphic 
distribution $W$ is integrable. Let $C \subset M$ be a  
standard minimal rational curve corresponding to a sufficiently
generic point of ${\Cal R}$.
Let ${\Cal W}_{C}$ be the closure of the leaf of $W$ containing $C$
and let ${\Cal W}_{1}$
be the closure of a leaf distinct from ${\Cal W}_{C}$.  Then, $C$ is disjoint
from ${\Cal S}(W)$ and $\Cal W_1$.
\endproclaim

\spb
\demo{Proof} Clearly, ${\Cal W}_{C} \cap {\Cal W}_{1}$ 
is contained in the singular locus ${\Cal S}(W)$.
So it is enough to show that a deformation of a 
standard minimal rational curve $C$ is disjoint from ${\Cal S}(W)$. 
This follows from a general result that a rational curve
with semipositive normal bundle can be deformed to avoid any codimension-2 set
(e.g. [Ko], Proposition 5.2.8).
We will give a proof for the reader's convenience.

\spb
For notational simplicity we assume that $C$ is embedded.
Let $[{\Cal O}(1)]^{p} \oplus {\Cal O}^{q}$ be the 
decomposition of the normal bundle $N_{C}$ of $C$ in $M$.
Choose a point $P$ on the curve away from ${\Cal S}(W)$. 
Choose sections $\sigma_{i}$, $1 \leq i \leq p$,
of the normal bundle
corresponding to independent sections of $[{\Cal O}(1)]^{p}$ vanishing at
$P$ and sections $\sigma_{i}$, $p+1 \leq i \leq n-1$ which generate the trivial
factor ${\Cal O}^{q}$ of the normal bunle. These sections $\sigma_{i}$, 
$1 \leq i \leq n-1$ are pointwise linearly independent outside $P$.
Consider an $(n-1)$-dimensional deformation of the curve corresponding to the
linear span of $\sigma_{i}$'s inside $H^{0}(C, N_{C})$, the tangent space of
deformations of $C$. Suppose all members of this $(n-1)$-dimensional 
family intersects ${\Cal S}(W)$. Since ${\Cal S}(W)$ has codimension 
$\geq 2$ in $M$,  this means that we have a 1-dimensional
subfamily intersecting ${\Cal S}(W)$ at a fixed point $Q \in {\Cal S}(W)$.
In particular, in the linear span of $\sigma_{i}$'s, 
there exists a section of the normal bundle vanishing at $Q \neq P$. 
But this is impossible because $\sigma_{i}$'s are pointwise independent
outside $P$. Thus, some member of this deformation must be disjoint
from ${\Cal S}(W)$. Proposition 12 follows.
\enddemo

\spb
We are now ready to deduce a topological obstruction to the integrability
of $W$ when cones are linearly degenerate.

\spb
\proclaim{Proposition 13} Suppose $b_{2}(M) = 1$ and the distribution $W$ is
non-trivial, i.e. ${\Cal C}_{x} \subset {\Bbb P}T_{x}(M)$ is linearly
degenerate for a generic $x \in M$. Then, $W$ cannot be integrable.
\endproclaim

\spb
\demo{Proof} Suppose $W$ is integrable. 
Let ${\Cal D}$ be an irreducible subvariety in the Chow space
of $M$ whose generic point corresponds to the closure of a generic leaf of $W$.
A generic hypersurface in ${\Cal D}$ corresponds to a hypersurface
$H$ in $M$ which is invariant under $W$, i.e. it is the union of leaves of
$W$. Let ${\Cal W}_{1}$ be the closure of a leaf of $W$. 
Then either ${\Cal W}_{1}  \subset  H$ or
${\Cal W}_{1} \cap H \subset {\Cal S}(W)$. From the previous proposition,
a generic standard minimal rational curve is disjoint
from $H$.  But from $b_{2}(M)=1$, $H$ is ample, a contradiction. 
\enddemo

\spb\noindent
{\smc Remarks}

\spb\noindent
With simple modifications the discussion in \S 4 applies even
if we drop the assumption that generic members of $\Cal K$ represent
standard minimal rational curves.

\bpb\noindent
{\bf \S 5 Cones in the central fiber}

\mpb\noindent
(5.1)\quad Returning to the situation of Theorem 1, we proceed to
examine the meromorphic distribution $W$ on the central fiber $X_0$
defined at generic points by cones $\Cal C_x\subset\Bbb PT_x(X_0)$. 
To start with, consider the irreducible Hermitian symmetric space $S$
with a base point $o\in S$. Let $\Cal C_o\subset\Bbb PT_o$ be the
standard cone. As in \S 4, let $\Cal T_o\subset\Bbb P\Lambda^2 T_o$
be the variety of tangential lines to $\Cal C_o$. We denote by
$\tilde{\Cal C}_o\subset T_o-\{0\}$ (resp. $\tilde{\Cal T}_o\subset
\Lambda^2 T_o-\{0\}$) the preimage of $\Cal C_o$ (resp. $\Cal T_o$)
under the canonical projection. We are going to prove

\spb
\proclaim{Proposition 14} $\Cal T_o \subset \Bbb P \Lambda^{2} T_{o}$ 
is linearly non-degenerate.
\endproclaim

\spb
For the proof we consider the principal bundle $\Cal G\subset\Cal F(S)$
defining  the standard $K^{\Bbb C}$-structure on $S$. It is based on a faithful
irreducible representation $\mu: K^{\Bbb C} \subset GL(V)$.
The cone ${\Cal C}_{o}$ can be identified
with the $K^{\Bbb C}$-orbit of projectivizations of 
highest weight vectors in ${\Bbb P}V$, via any choice of an element of
$\Cal G\subset\text{Isom}(V,T_0)$.

\spb
We recall Grothendieck's classification of $G$-principal bundles
over ${\Bbb P}^{1}$, where $G$ is any connected complex reductive Lie group.
Choose a maximal algebraic torus  $H \subset G$.
Let ${\Cal O}(1)^{*}$ be the ${\Bbb C}^{*}$-principal
bundle on ${\Bbb P}^{1}$, which is just the complement of the zero section
of ${\Cal O}(1)$. 

\spb
\proclaim{Proposition 15 (Grothendieck [Gro])}
Let ${\Cal G}$ be a principal $G$-bundle on ${\Bbb P}^{1}$. Then, there exists an algebraic 
one-parameter subgroup $\rho: {\Bbb C}^{*} \rightarrow H$ such that ${\Cal G}$
is equivalent to the $G$-bundle associated to ${\Cal O}(1)^{*}$ via the action
$\rho$. Furthermore, let ${\Cal V}$ be a vector bundle associated to ${\Cal G}$
via a representation $\mu : G \rightarrow GL(V)$. Then ${\Cal V}$ splits
as the direct sum of line bundles ${\Cal O}(<\mu_{i}, \rho>)$, where
$\mu_{i}: H \rightarrow {\Bbb C}^{*}$ are the weights of $\mu$ and $<\mu_{i},
\rho>$ denotes the integral exponent of the homomorphism
$\mu_{i} \circ \rho: {\Bbb C}^{*} \rightarrow {\Bbb C}^{*}$.
\endproclaim

\spb
\demo{Proof of Proposition 14}
We apply Grothendieck's result 
to the case $G= K^{{\Bbb C}}$.  For any $\alpha \in 
{\Cal C}_{x}$, we have a minimal rational curve $C$ passing through $o \in S$, 
tangential to $\alpha$ with $T(S)|_{C} = {\Cal O}(2) \oplus [{\Cal O}(1)]^{p}
\oplus {\Cal O}^{q}$ for some positive integers $p,q$.
Restricting ${\Cal G}$ to $C$, we obtain a $K^{{\Bbb C}}$-principal bundle on
$C$. From  Grothendieck's classification above,
 we have a one-parameter subgroup $\rho:
{\Bbb C}^{*} \rightarrow K^{{\Bbb C}}$ defining
 ${\Cal G}|_{C}$. From the splitting
type of $T(S)|_{C}$, $\rho({\Bbb C}^{*})$ acts on $V$ with three exponents
$2, 1, 0$ and gives rise to a decomposition $V = {\Bbb C} \alpha \oplus
{\Cal H}_{\alpha} \oplus {\Cal N}_{\alpha}$. This action preserves the
cone ${\Cal C}_{o} \subset {\Bbb P} V$.

\spb
Taking the inverse of $\rho$ and multiplying by a scalar representation, 
we get a ${\Bbb C}^{*}$-representation on $V$ preserving ${\Cal C}_{o}$ which
fixes ${\Bbb C} \alpha$, acts as $t$ on ${\Cal H}_{\alpha}$ and acts as
$t^{2}$ on ${\Cal N}_{\alpha}$, $t \in {\Bbb C}^{*}$. Choose a generic point
$\alpha + \xi + \zeta$ on $\tilde {\Cal C}_{o}$, 
with $\xi \in {\Cal H}_{\alpha}$
and $\zeta \in {\Cal N}_{\alpha}$. The orbit of the ${\Bbb C}^{*}$-action 
is $\alpha + t \xi + t^{2} \zeta$. At $t = t_{0}$, we further
consider the curve $\alpha + e^{s} t_{0} \xi + e^{2s} t_{0}^{2} \zeta$, $s \in
{\Bbb C}$.  Taking derivative with respect to $s$ at $s=0$, we get 
the tangent vector $t_{0}\xi + 2 t_{0}^{2} \zeta$ to $\tilde{\Cal C}_{o}$
at the point $\alpha + t_{0} \xi + t_{0}^{2} \zeta$. The corresponding
element of $\tilde{\Cal T}_{o}$ is 
$$
(\alpha + t_{0} \xi + t_{0}^{2} \zeta) \wedge (t_{0} \xi + 2 t_{0}^{2} \zeta)
= t_{0} \alpha \wedge \xi + 2 t_{0}^{2} \alpha \wedge \zeta +t_{0}^{3}
\xi \wedge \zeta\ .
$$
Thus, the linear span of ${\Cal T}_{o}$ contains $\alpha \wedge \xi, \alpha
\wedge \zeta$, and $\xi \wedge \zeta$ for any generic $\alpha \in {\Cal C}_{o}$
and $\alpha + \xi + \zeta \in {\Cal C}_{o}$. As we vary $\xi$ on 
${\Cal H}_{\alpha}$, the corresponding vectors 
$\zeta \in {\Cal N}_{\alpha}$ span
${\Cal N}_{\alpha}$. Otherwise, ${\Cal C}_{o}$ will be contained in the linear
subspace ${\Bbb C} \alpha \oplus {\Cal H}_{\alpha} \oplus {\Cal N}'_{\alpha}$
for some proper subspace ${\Cal N}'_{\alpha} \subset {\Cal N}_{\alpha}$,
contradicting the linear non-degeneracy of ${\Cal C}_{o}$. It follows that
${\Cal T}_{o}$ is linearly non-degenerate.
\enddemo

\spb\noindent
{\smc Remarks}

\spb\noindent
The above proposition can be proved by checking case by case using
representation theory. Let $\frak l= [\frak k^{\Bbb C}, \frak k^{\Bbb C}]$
be the maximal semisimple subalgebra. 
The representation of $\frak l$ on $ T_{o}$  can be identified
from the table in Section 2. In Table 5 of 
Onishchik-Vinberg [OV], the irreducible decomposition
of the second exterior power of these representations can be found.
We will list them below. The notations are as in [OV],
Table 1 (pp.~292-294) and Table 5 (pp.~299-305). $\pi_{i}$ denotes the $i$-th 
fundamental weight of the simple Lie algebra and $R(\lambda)$ denotes the
irreducible representation with the highest weight $\lambda$. 
For $\frak{sl}(p)\times\frak{sl}(q)$, $R(\lambda) \otimes R(\mu)$ 
denotes the tensor product
of the $\frak{sl}(p)$-representation $R(\lambda)$ 
and the $\frak{sl}(q)$-representation $R(\mu)$.

\newpage
$$
\vbox{\offinterlineskip
\halign to \displaywidth{\tabskip=0pt plus1fil
\vrule#
&\strut\hfil$#$\hfil
&\vrule#
&\hfil$#$\hfil
&\vrule#
&\hfil$#$\hfil
&\vrule#
&\hfil$#$\hfil
&\vrule#
\tabskip=0pt\cr
\multispan{9}\hrulefill\cr
height3pt&\omit&&&&&&&\cr
&\text{Type}&&l&&T_o&&\Lambda^2 T_o&\cr
height3pt&\omit&&&&&&&\cr
\multispan{9}\hrulefill\cr
height3pt&\omit&&&&&&&\cr
&\text{I}&&\frak{sl}(p)\times\frak{sl}(q)&&R(\pi_1)\otimes R(\pi_1)&&R(\pi_2)\otimes R(2\pi_1)\oplus R(2\pi_1)\otimes R(\pi_2)&\cr
height3pt&\omit&&&&&&&\cr
&\text{II}&&\frak{sl}(n)&&R(\pi_2)&&R(\pi_1+\pi_3)&\cr
height3pt&\omit&&&&&&&\cr
&\text{III}&&\frak{sl}(n)&&R(2\pi_1)&&R(2\pi_1+\pi_2)&\cr
height3pt&\omit&&&&&&&\cr
&\text{IV}&&\frak{so}(n)&&R(\pi_1)&&R(\pi_2)&\cr
height3pt&\omit&&&&&&&\cr
&\text{V}&&\frak{so}(10)&&R(\pi_5)&&R(\pi_3)&\cr
height3pt&\omit&&&&&&&\cr
&\text{VI}&&\frak e_6&&R(\pi_1)&&R(\pi_2)&\cr
height3pt&\omit&&&&&&&\cr
\multispan{9}\hrulefill\cr
}}$$

\spb
With the exception of Type I, 
$\Lambda^{2} T_{o}$ is irreducible.
Proposition 14 is obvious in these cases, 
because the linear span of $\Cal T_{o}$ is an invariant subspace.

\spb
For Type I, we can check it in the following way. Let $v \in T_{o}$ be the
highest weight vector with weight $\pi_{1} \otimes \pi_{1}$.
There are two weights maximal among the remaining weights, $\pi_{1} \otimes 
(\pi_{2} - \pi_{1})$ and $(\pi_{2} - \pi_{1}) \otimes \pi_{2}$. Corresponding
weight vectors are $l \cdot v$ and $l' \cdot v$ respectively, where $l$ (resp.
$l'$ ) denotes an eigenvector of $\frak{sl}(p)$ (resp. $\frak{sl}(q)$) 
corresponding to the root $\alpha_{1}$. 
Hence the highest weight
vectors of $\Lambda^{2} T_{o}$ are 
$v \wedge (l \cdot v)$ and $ v \wedge (l' \cdot v)$. But 
$v \wedge(l \cdot v)$ and $ v \wedge (l' \cdot v)$ correspond to
tangential lines to the cone at $v$. Hence both highest
weight vectors are in ${\Cal T}_{o}$, as desired. 

\spb
To study the meromorphic distribution $W$ on $X_0$ we prove 

\spb
\proclaim{Proposition 16} 
Let $W_{o} \subset T_{o}$ be a subspace and ${\Bbb P}T_{o}
\rightarrow {\Bbb P} W_{o}$ be a projection defined by 
choosing a complementary subspace. Suppose 
the strict image of ${\Cal C}_{o}$ is linearly
non-degenerate in ${\Bbb P}W_{o}$. Then, the tangential lines to 
the smooth part of the strict image generates ${\Bbb P} \Lambda^{2} W_{o}$.
\endproclaim

\spb
\demo{Proof} The projection $T_{o} \rightarrow W_{o}$ induces a surjective
map $\Lambda^{2} T_{o} \rightarrow \Lambda^{2} W_{o}$. Proposition 16
is thus a direct consequence of Proposition 14. 
\enddemo

\mpb\noindent
(5.2)\quad We are ready to finalize the proof of Theorem 1, for which we
may assume that $S$ is of $\text{rank}\ge 2$. We apply results of \S 4
to $M=X_0$. The family $\Cal K$ is chosen as the family $\Cal K$ of (2.2).
Then, the irreducibility assumption made in \S 4 holds. 
By [(2.4), Proposition 6] a generic member of $\Cal K$ represents
a standard minimal rational curve.
$W_x$, $\Cal C_x$ and $\Cal T_x$ are defined for a generic point 
$x \in X_{0}$. Let $\sigma: \Delta \rightarrow  X$ be a section of $\pi$, 
so that $\sigma(0) = x$. From \S 3, we know that the family ${\Cal M}_{\sigma(t)}$ is 
a trivial family of standard cones ${\Cal C}_{o}$. We have the family of
rational maps $\Theta_{\sigma(t)}: {\Cal M}_{\sigma(t)} \rightarrow {\Bbb P}
T_{\sigma(t)}(X_{t})$. Thus, $\Theta_{x}: {\Cal M}_{x}\cong {\Cal C}_{o}
\rightarrow {\Bbb P}T_{x} (X_{0})$ is a (generically finite)
rational map, defined by a subsystem of the
complete linear system defining the standard embedding of ${\Cal C}_{o}$.
By applying [(5.1), Proposition 16], 
${\Cal T}_{x}$ is linearly non-degenerate in
${\Bbb P} \Lambda^{2} W_{x}$. By [(4.2), Proposition 9],
$W$ is integrable. Since $b_2(X_0)=1$, by [(4.3), Proposition 13]
$W$ must be the trivial distribution $T_{X_0}$. Thus,
$\Cal C_x\subset\Bbb PT_x(X_0)$ is linearly non-degenerate at a
generic point $x$ and is the standard embedding given by the complete
linear system. By the Hartogs' extension result
[(1.2), Proposition 3] there exists on $X$ a holomorphic bundle of
cones $\Cal C_x\subset\Bbb PT_x$, defined 
everywhere on $X$, such that the inclusion is isomorphic to the
embedding of the standard cone $\Cal C_o\subset\Bbb PT_o(S)$ on the
model space $S$. Thus, we have obtained a holomorphic $K^{\Bbb C}$-structure
on $X_0$, which is integrable by the closedness of the integrability
condition (cf. [(1.4), Proposition 2]). As $X_0$ is simply-connected,
it follows from the result of Ochiai [(1.4), Proposition 3] that $X_0$ is
biholomorphic to $S$. The proof of Theorem 1 is completed.

\spb
Finally, for Theorem 1$'$ we note that $\rho$ is projective because
$\text{Pic}(S)\cong\Bbb Z$. Then,
Theorem 1$'$ follows from Theorem 1 by slicing by algebraic curves
$Z_0$ on the base space $Z$, noting that if the fiber $X_y$ is biregular
to $S$ at some point $y$ on the curve $Z_0$, then $X_z$ is biregular to
$S$ for all $z\in Z_0$ except for at worst a finite number of base points
$z$, since $S$ is infinitesimally rigid. 

\spb
For completeness we explain here why the apparently stronger Theorem 1
is {\it a priori\/} equivalent to Theorem 1$'$. It suffices to show that 
modulo at worst a lifting by a ramified  covering $\triangle\to\triangle$,
the regular family $\pi:X\to\triangle$ can be completed to a regular
family $\pi':X'\to\Bbb P^1$ over the Riemann sphere $\Bbb P^1$ such that
the total space $X'$ is projective-algebraic. The extension can in fact be
obtained by gluing the regular family $\pi:X\to\triangle$ with the
trivial family on $D=\Bbb P^1-\overline{\triangle(\frac 12)}$, provided
that the regular family $\pi:X\to\triangle$ is holomorphically trivial
on the annulus $A=\triangle\cap D=\{z:\frac 12<|z|<1\}$. The structure group
of $\pi|_A:X|_A\to A$, $X|_A=\pi^{-1}(A)$, is the group $\text{\it Aut}(S)$
of automorphisms of $S$. The identity component $\text{\it Aut}_0(S)$ is a 
complex Lie group and $\text{\it Aut}(S)/\text{\it Aut}_0(S)$ is finite.
Thus, lifting by at worst a finite ramified covering 
$\varphi:\triangle\to\triangle$ defined by $\varphi(z)=z^p$,
we may assume that the structure group is the connected complex 
Lie group $\text{\it Aut}_0(S)$. It then follows that $\pi|_A:X|_A\to A$ is
topologically trivial and hence holomorphically trivial by the 
Oka-Grauert Principle (Grauert [G], Cartan [C]) 
for holomorphic principal bundles,
allowing us to define an extension $\pi':X'\to\Bbb P^1$ by gluing.
Finally we observe that $X'$ is projective-algebraic. Let $\Lambda$ be the
determinant bundle of the relative tangent bundle $T$ on $X'$. As all
fibers are K\"ahler, $\Lambda|_{X'_t}$ is ample for any $t\in\Bbb P^1$.
Let $H$ be the hyperplane section line bundle on $\Bbb P^1$.
Then, $\Lambda\otimes(\pi^{\prime *}H)^k$ is ample for $k$ sufficiently
large, as desired.

\spb\noindent
{\smc Acknowledgements.}\quad Part of the work was done while the
first author visited The University of Hong Kong in the summer of
1995. He would like to thank The University of Hong Kong and University
of Notre Dame for the support for the visit. Both authors learned
about the problem from Professor Y.-T. Siu, to whom they would like to
express their thanks. They also like to thank Sai-Kee Yeung for his
interest and for providing the reference [OV], and Professor S. Kobayashi
for clarifying the history regarding $G$-structures.

\newpage

\bpb\noindent
{\bf References}

\spb
\item{[Bo]} Bott, R.\quad
Homogeneous vector bundles, 
Annals of Math. {\bf 66} (1957), 203-248.

\item{[Br]} Brieskorn, E.\quad
Ein Satz \"uber die komplexen Quadriken,
Math. Ann. {\bf 155} (1964), 184-193.

\item{[C]} Cartan, H.\quad
Espaces fibr\'es analytiques, Sym. Int. de Top. Alg.,
Univ. Nac. de Mexico 1958, pp.97-121.

\item{[Gra]} Grauert, H.\quad
Holomorphe Funktionen mit Werten in Komplexen Lieschen Gruppen,
Math. Ann. {\bf 133} (1957), 450-472.

\item{[Gro]} Grothendieck, A.\quad
Sur la classification des fibr\'es holomorphes sur la sph\`ere de Riemann.
Amer. J. Math. {\bf 79} (1957), 121-138.

\item{[HK]} Hirzebruch, F. and Kodaira, K.\quad
On the complex projective spaces,
J. Math. Pures Appl. {\bf 36} (1957), 201-216.

\item{[Ho]} Holmann, H.\quad
Komplexe R\"aume mit komplexen Transformationsgruppen,
Math. Ann. {\bf 150} (1963), 327-360.

\item{[H1]} Hwang, J.-M.\quad
Nondeformability of the complex hyperquadric, 
Invent. math. {\bf 120} (1995), 317-338.

\item{[H2]} Hwang, J.-M.\quad
Characterization of the complex projective space by holomorphic vector fields,
to appear in Math. Zeit.

\item{[KS]} Kodaira, K. and Spencer, D.\quad
On deformations of complex structures II,
Ann. Math. {\bf 67} (1958), 403-466.

\item{[Ko]} Koll\'ar, J. \quad
Flips, flops, minimal models, etc., 
Surveys in Differential 
\linebreak
Geometry {\bf 1} (1991), 113-199.

\item{[Ma]} Mabuchi, T.\quad
K\"ahlerian rigidity of irreducible Hermitian symmetric spaces. Preprint 1993.

\item{[MM]} Matsushima, Y. and Morimoto, A.\quad
Sur certains espaces fibr\'es holomorphes sur une vari\'et\'e de Stein,
Bull. Soc. Math. France {\bf 88} (1960), 137-155.

\item{[Mk1]} Mok, N.\quad
The uniformization theorem for compact K\"ahler manifolds of nonnegative
holomorphic bisectional curvature, 
J. Diff. Geom. {\bf 27} (1988), 179-214.

\item{[Mk2]} Mok, N.\quad
{\it Metric Rigidity Theorems on Hermitian Locally Symmetric Manifolds\/}
(Ser. Pure Math., Vol.6), World Scientific, 
Singapore-New Jersey-London-Hong Kong 1989.

\item{[Mr]} Mori, S.\quad
Projective manifolds with ample tangent bundles, 
Ann. Math. {\bf 110} (1979), 593-606.

\item{[NT]} Nakagawa, H. and Takagi, R.\quad
On locally symmetric K\"ahler manifolds in a complex projective space,
J. Math. Soc. Japan {\bf 28} (1976), 638-667.

\item{[Oc]} Ochiai, T. \quad
Geometry associated with semisimple flat homogeneous spaces. 
Trans. A.M.S. {\bf 152} (1970), 159-193.

\item{[OV]} Onishchik, A. and Vinberg, E.\quad
{\it Lie Groups and Algebraic Groups\/},
Springer Verlag, Berlin-Heidelberg-New York, 1990.

\item{[R]} Ros, A.\quad
A characterization of seven compact K\"ahler submanifolds 
by holomorphic pinching, 
Ann. Math. {\bf 121} (1985), 377-382.

\item{[S1]} Siu, Y.-T.\quad
Nondeformability of the complex projective space,
J. reine angew. Math. {\bf 399} (1989), 208-219.
Errata {\bf 431} (1992), 65-74.

\item{[S2]} Siu, Y.-T.\quad
Uniformization in Several Complex Variables,
in {\it Contemporary Geometry\/}, J.Q. Zhong Memorial Volume,
ed. H.-H. Wu, Plenum Press New York-London, 1991.

\item{[St]} Sternberg, S. \quad
{\it Lectures on Differential Geometry\/}, 
Prentice-Hall, Englewood Cliffs, New Jersey, 1964.

\item{[W]} Wolf, J.A.\quad
Fine structures of Hermitian symmetric spaces,
in {\it Geometry of Symmetric Spaces\/}, ed. Boothby-Weiss, Marcel-Dekker,
New York 1972, pp.271-357.

\address
Jun-Muk Hwang,
University of Notre Dame, Notre Dame, IN 46556
\endaddress

\curraddr
MSRI, 1000 Centennial Drive, Berkeley, CA 94720 
\endcurraddr

\email
jhwang\@msri.org
\endemail

\address
Ngaiming Mok,
The University of Hong Kong, Pokfulam, Hong Kong
\endaddress

\email
nmok\@hkucc.hku.hk
\endemail

\enddocument